\documentclass[11pt, reqno]{amsart}

\usepackage{amssymb}
\usepackage{hyperref} 
\usepackage{mathrsfs,bbold}
\usepackage{amsmath}
\usepackage{amsthm}
\usepackage{enumerate}
\usepackage{dsfont}
\usepackage{color}
\usepackage[a4paper]{geometry}
\usepackage[utf8]{inputenc}
\usepackage{stmaryrd}
\usepackage{titlesec}
\usepackage{mathabx}
\usepackage{appendix}
\usepackage{todonotes, fancyhdr}

\usepackage[all]{xy}
\let\oldtocsection=\tocsection
\let\oldtocsubsection=\tocsubsection
\let\oldtocsubsubsection=\tocsubsubsection

\renewcommand{\tocsection}[2]{\hspace{0em}\oldtocsection{#1}{#2}}
\renewcommand{\tocsubsection}[2]{\hspace{1em}\oldtocsubsection{#1}{#2}}
\renewcommand{\tocsubsubsection}[2]{\hspace{2em}\oldtocsubsubsection{#1}{#2}}

\numberwithin{theorem}{section}
\numberwithin{equation}{section}

\newcommand{\ssk}{\smallskip}

\renewcommand{\epsilon}{\varepsilon}
\newcommand{\eps}{\epsilon}

\newcommand\bbE{\mathbb{E}}

\newcommand\bbN{\mathbb{N}}
\newcommand\bbR{\mathbb{R}}

\newcommand\bfR{\textsf{\textbf{R}}}

\newcommand{\bfX}{\textsf{\textbf{X}}}

\newenvironment{Dem}[1][\unskip]{%
    \begin{list}{\hspace{0.5cm}{\sf \textbf{Proof #1 --}}}{%
        \setlength{\topsep}{0pt}%
        \setlength{\leftmargin}{0pt}%
        \setlength{\rightmargin}{0pt}%
        \setlength{\listparindent}{0pt}%
        \setlength{\itemindent}{0pt}%
        \setlength{\parsep}{0pt}%
        \addtolength{\leftmargin}{20pt}%
        \addtolength{\rightmargin}{0pt}%
    } \item }{\hfill $\rhd$\end{list}\smallskip}

\titleformat{\section}[block]
{\filcenter\normalfont\sffamily\bfseries\Large}
{{\hspace{-0.7cm}}\thesection \hspace{0.2em} --\vspace{0.3cm}}{0.5em}{}

\titleformat{\subsection}[block]
{\filcenter\normalfont\sffamily\bfseries\large}  						  
{\hspace{-0.7cm}\thesubsection \hspace{0.5em}--\vspace{0.3cm}}{.5em}{}  
\titlespacing{\subsection}{-0pc}{1.5ex plus .1ex minus .2ex}{0pc}

\titleformat{\subsubsection}[block]
{\normalfont\sffamily\bfseries}{\hspace{-1cm}\thesubsubsection}{.5em}{}
\titlespacing{\subsubsection}{15pc}{1ex plus .1ex minus .2ex}{1pc}

\def\XXint#1#2#3{{\setbox0=\hbox{$#1{#2#3}{\int}$}
     \vcenter{\hbox{$#2#3$}}\kern-.5\wd0}}

\numberwithin{subsection}{section}
\numberwithin{subsubsection}{subsection}

\newtheoremstyle{mystyle}
{3pt}               
{3pt}               
{\it }                      
{}                      
{\sffamily\bfseries}             
{}                      
{0.5em}                 
{\llap{#2. }#1{$\;$ --}}

\theoremstyle{mystyle}

\newtheorem{thm}{Theorem}
\newtheorem*{thm*}{Theorem}

\newtheorem{cor}[thm]{\hspace{-0.15cm}  {Corollary} }
\newtheorem{lem}[thm]{\hspace{-0.2cm}  {Lemma} }
\newtheorem{prop}[thm]{\hspace{-0.2cm} {Proposition}}
\newtheorem{defn}[thm]{ \hspace{-0.3cm} {Definition}}
\newtheorem*{defn*} {Definition}
\newtheorem*{prop*} {Proposition}
\newtheorem*{lem*} {Lemma}
\newtheorem*{cor*} {Corollary}
\newtheorem{rem}[thm]{\hspace{-0.15cm} {Remark}}
\newtheorem{assum}[thm]{Assumption}

\newtheoremstyle{mystyle2}
{3pt}               
{3pt}               
{\it }                      
{}                      
{\sffamily\bfseries}             
{}                      
{0.5em}                 
{\llap{#2 }#1{\hspace{0.2cm}--}}
\theoremstyle{mystyle2}

\newtheorem*{definition*}{Definition}
\newtheorem*{theorem*}{Theorem}
\numberwithin{equation}{section} 
\newtheorem*{Remark*}{Remark}

\newcommand\restr[2]{{
  \left.\kern-\nulldelimiterspace 
  #1 
  \vphantom{\big|} 
  \right|_{#2} 
  }}

\newcommand{\bX}{\mathbf{X}}

\newcommand{\RR}{\mathbb{R}}

\newcommand{\fI}{\mathfrak{I}}

\newcommand{\dd}{\mathrm{d}}

\begin{document}

\vspace*{3ex minus 1ex}
\begin{center}
{\Huge\sffamily{Non-explosion criteria for rough differential equations driven by unbounded vector fields   \vspace{0.5cm}}}
\end{center}
\vskip 5ex minus 1ex

\begin{center}
{\sf I. BAILLEUL\footnote{I.Bailleul thanks the U.B.O. for their hospitality, part of this work was written there.} and R. CATELLIER\footnote{R. Catellier acknoledges the support of the Lebesgue Centre.}}
\end{center}

\vspace{1cm}

\begin{center}
\begin{minipage}{0.8\textwidth}
\renewcommand\baselinestretch{0.7} \scriptsize \textbf{\textsf{\noindent Abstract.}} 
We give in this note a simple treatment of the non-explosion problem for rough differential equations driven by unbounded vector fields and weak geometric rough paths of arbitrary roughness.   
\end{minipage}
\end{center}

\vspace{1cm}

\section{Introduction}

Although rough paths theory has now been explored for twenty years, a few elementary questions are still begging for a definite answer. We consider the existence problem for the local time and occupation time of solutions to rough differential equations as the main open problem, in relation with reflection problems. At a more fundamental level, the question of global in time existence of solutions of a rough differential equation
\begin{equation}
\label{EqRDE}
dz_t = \textrm{F}(z_t)\,d\bfX_t,
\end{equation}
under relaxed boundedness assumptions on the vector fields $\textrm{F}= (V_1,\dots,V_\ell)$ has not been clarified so far. Given a weak geometric $p$-rough path $\bfX$ defined on some time interval $[0,T]$, the preceding equation is known to have a solution defined on the whole of $[0,T]$ if the driving vector fields $V_i$ are $C^\gamma_b$, for some regularity exponent $\gamma>p$; see for instance T. Lyons' seminal paper \cite{Lyons98} or the lecture notes \cite{LyonsStFlour}. One would ideally like to relax these boundedness assumptions to some linear growth assumption, but the following elementary counter-examples of Gubinelli and Lejay \cite{LejayGlobal} shows that this is not sufficient. Consider the dynamics \eqref{EqRDE} on $\RR^2$, with $\textrm{F} = (V_1,V_2)$, and vector fields $V_1(x,y) = (x\sin(y),x)$ and $V_2(x,y) = 0$, driven by the non-geometric pure area rough path ${\bf X}_t = 1+t(1\otimes 1)$. Writing $z_t = (x_t,y_t)$, one sees that $z$ is actually the solution of the ordinatry differential equation 
$$
\dot z_t = (\dot x_t,\dot y_t) = \big(x_t\sin(y_t)^2 + x_t^2\cos(y_t),\, x_t\sin(y_t)\big).
$$
The solution started from an initial condition of the form $(a,0)$, with $a$ positive, has constant null $y$-component and has an exploding $x$-component since $\dot x_t = x_t^2$.

\medskip

The non-explosion problem was explored in a number of works for differential equations driven by $p$-rough paths, for $2\leq p<3$, especially in the works of Davie \cite{Davie} and Lejay \cite{Lejay, LejayGlobal}. Davie provides essentially the sharpest result in the regime $2\leq p<3$.   \vspace{0.1cm}

\begin{itemize}
   \item To make it simple, assume F is $C^3$ and has linear growth: $\big|\textrm{F}(x)\big| \lesssim |x|$. Theorem 6.1 (a) in \cite{Davie} provides a non-explosion criterion in terms of the growth rate of $D^2\textrm{F}$
$$
\big|D^2F(x)\big| \leq h(|x|).
$$
There is no explosion if $h(r)\lesssim \frac{1}{r}$, and
$$
\int^\infty \left(\frac{r^{\gamma-2}}{h(r)}\right)^{\frac{p-1}{\gamma-1}}\,\frac{dr}{r^p} = \infty.
$$
Davie's criterion is shown to be sharp in the class of all $p$-rough paths, $2\leq p<3$, with an example of a rough differential equation where explosion can happen for some appropriate choice of a non-weak geometric rough path in case the criterion is not satisfied -- see Section 6 in \cite{Davie}. The limit case for Davie's criterion is $h(r) = \frac{O(1)}{r}$. We essentially recover that bound.   \vspace{0.15cm}

   \item Lejay \cite{Lejay} works with Banach space valued weak geometric $p$-rough paths, with $2\leq p<3$. In the setting where the vector fields $V_i$ are $C^3$ with bounded derivates and are required to have growth rate $\big|V_i(x)\big| \lesssim g(|x|)$, he shows non-explosion of solutions to equation \eqref{EqRDE} under the condition that $\sum_k \frac{1}{g(k)^p}$ diverges. The limit case is $g(r)\simeq r^\frac{1}{p}$.   \vspace{0.15cm}
   
   	\item The analysis of Friz and Victoir \cite{FrizVictoirBook}, Exercice 10.56, gives a criterion comparable to ours, with an erronous proof. They use a pattern of proof that is implemented in a linear setting and cannot work in a nonlinear framework as it bears heavily on a scaling argument -- see the proof of Theorem 10.53. One can see part of the present work as a correct or alternative proof of their statement.
\end{itemize}

\medskip

We identify in the sequel a vector field $V$ on $\RR^d$ with the first order differential operator $f\mapsto (Df)(V)$. For a tuple $I = (i_1,\dots,i_k)\in\{1,\dots,\ell\}^k$, and vector fields $V_1,\ldots, V_\ell$, we define the differential operators
$$
V_I := V_{i_1}\cdots V_{i_k}, \quad \textrm{and}\quad V_{[I]} := \Big[V_{i_1},\dots, [V_{i_{k-1}},V_{i_k}]\big]\Big],
$$
under proper regularity assumptions on the $V_i$. (Note that the operator $V_{[I]}$ is actually of order one, so $V_{[I]}$ is a vector field.) The local increment $z_t-z_s$ of a solution $z$ to the rough differential equation \eqref{EqRDE} is known to be well-approximated by the time 1 value of the ordinary differential equation
\begin{equation}
\label{EqApproximateODE}
y'_r = \sum_{k=1}^{[p]} \sum_{I\in\left\{1,\cdots,\ell\right\}^k} \Lambda^{k,I}_{ts} V_{[I]}\big(s,y_r(x)\big),
\end{equation}
where ${\bf \Lambda}_{ts} := \log {\bf X}_{ts}$, and $0\leq r\leq 1$ -- see \cite{BailleulRMI} or \cite{BoutaibGyurkoLyonsYang} for instance. The following simplified version of our main result, Theorem \ref{thm:main}, actually gives a non-explosion result in terms of growth assumptions on the vector fields $V_{[I]}$ that appear in the approximate dynamics \eqref{EqApproximateODE}. Pick an arbitrary $p>1$ and a weak geometric $p$-rough path $\bf X$.

\medskip

\begin{thm}
\label{ThmPreliminary}
There is no explosion for the solutions of the rough differential equation \eqref{EqRDE} is the functions $V_{i_1}\cdots V_{i_n}\textrm{\emph{Id}}$ are $C^2$ with bounded derivatives, for any $1\leq n\leq [p]$ and any tuple $(i_1,\ldots,i_n)\in \llbracket 1,\ell\rrbracket^n$.
\end{thm}

\medskip

Theorem \ref{thm:main} is sharper than that statement as it involves the vector fields $V_{[I]}$ -- recall Example 3 of \cite{LejayGlobal}. In the case where $2\leq p<3$, our non-explosion criterion becomes
$$
\big|D^2_x\textrm{F}\big| \vee \big|D^3_x\textrm{F}\big| \lesssim \frac{1}{1+|x|},
$$
for a multiplicative implicit constant independent of $x\in\RR^d$. We mention here that we have been careful on the growth rate of the different quantities but that one can optimize the regularity assumptions that are made on the vector fields $V_i$ to get slightly sharper results. This explains the discrepancy between Davie's optimal criterion in the case $2\leq p<3$ and our result. These refinements are not needed for the applications \cite{BailleulCatellierMeanField}; we leave them to the reader. Note also here that one can replace $\RR^d$ by a Banach space and give versions of the statements involving infinite dimensional rough paths, to the price of using slightly different notations, such as in \cite{BailleulSeminaire}. There is no difference between the finite and the infinite dimensional settings for the explosion problem.

\bigskip

Our main result, Theorem \ref{thm:main}, holds for dynamics \eqref{EqRDE} with a drift and time-dependent vector fields. It is proved in Section \ref{SectionFlows} on the basis of some intermediate technical estimates whose proof is given in Section \ref{SectionProof}. Theorem \ref{thm:main} holds for H\"older $p$-rough paths. A similar statement holds for more general continuous rough paths, with finite $p$-variation, such as proved in Section \ref{SectionCorollaries} with other corollaries and extensions.

\bigskip

\noindent \textsf{\textbf{Notations.}} We gather here a number of notations that are used throughout the paper.   \vspace{0.1cm}

\begin{itemize}
      \item Given a positive finite time horizon $T$, we denote by $\Delta_T$ the simplex $\{(t,s)\in [0,T]^2\, : \, 0 \le s \le t \le T\}$.   \vspace{0.1cm}
   
   \item We refer the reader to Lyons' seminal article \cite{Lyons98} or any textbook or lectures notes on rough paths \cite{LyonsQian, LyonsStFlour, FrizVictoirBook, BaudoinBookEMS, LNRoughPaths} for the basics on rough paths theory and simply mention here that we work throughout with finite dimensional weak geometric H\"older $p$-rough paths ${\bf X} = 1\oplus X^1\oplus \cdots \oplus X^{[p]}$, with values in $\bigoplus_{i=0}^{[p]}({\bfR}^\ell)^{\otimes i}$ say, and norm
   $$
   \|\bX\| := \max_{1\leq i\leq [p]} \; \sup_{0\leq s < t \leq T} \frac{\big|X^i_{ts}\big|^{\frac1i}}{|t-s|^{\frac{1}{p}}}.
   $$
   Note that if $\mathbf{\Lambda} = \big(0\oplus \Lambda^1 \oplus \cdots \oplus \Lambda^{[p]}\big)$ is the logarithm of the rough path $\bX$, we have for all $0 \le s \le t \le T$, all $i\in\{1,\cdots,[p]\}$, 
\[
\big|\Lambda^i_{ts}\big| \lesssim_{i} \|\bX\|^i |t-s|^{\frac{i}{p}}
\]   

   \item Last, we use the notation $a\lesssim b$ to mean that $a$ is smaller than a constant times $b$, for some universal numerical constant.
\end{itemize}

\bigskip

\section{Solution flows to rough differential equations}
\label{SectionFlows}

Pick $\alpha\in [0,1]$. A finite dimensional-valued function $f$ defined on $\RR^d$ is said to have $\alpha$-growth if
$$ 
\sup_{x\in\RR^d}\frac{\big|f(x)\big|}{\big(1+|x|\big)^\alpha} < +\infty.
$$ 
Let $V_0$ and $V_1,\cdots,V_d : [0,T]\times \RR^d \to \RR^d$ be time-dependent vector fields on $\RR^d$. 

\medskip

\begin{assum}
\label{hyp:1}
\emph{\textsf{\textbf{Space regularity and growth.}}} For any $1\leq n\leq [p]$ and for any tuple $I\in \{1,\cdots ,\ell\}^n$,
\begin{enumerate}
   \item[$\bullet$] the vector fields $V_0(s,\cdot)$ and $V_{[I]}(s,\cdot)$ are Lipschitz  continuous with $\alpha$-growth, and their derivatives 	$DV_0(s,\cdot)$ and $DV_{[I]}$ are $C^1_b(B,B)$, uniformly in time,   \vspace{0.1cm}
	
   \item[$\bullet$] for all indices $1\leq k_1,\cdots,k_n \leq [p]$ with $\sum k_i \le [p]$, and all tuples $I_{k_i}\in\{1,\cdots,\ell\}^{k_i}$, the functions
   $$ 
   V_0(s,\cdot)V_{[I_{n-1}]}(s,\cdot)\cdots V_{[I_1]}(s,\cdot)\textrm{\emph{Id}} \quad \text{and} \quad V_{[I_{n}]}(s,\cdot)\cdots V_{[I_1]}(s,\cdot)\textrm{\emph{Id}}
   $$
   are $C^2_b$ with $\alpha$-growth, uniformly in time.
\end{enumerate}
\end{assum}

\medskip

One can trade in the above assumption some growth condition on the $V_i$ against some growth condition on its derivatives; this is the rationale for introducing the notion of $\alpha$-growth.

\medskip

\begin{assum}
\label{hyp:time}
\emph{\textsf{\textbf{Time regularity and growth.}}} There exists some regularity exponents $\kappa_1\ge\frac{1+[p] - p}{p}$ and $\kappa_2 \ge\frac{[p]}{p}$ with the following properties. 
\begin{itemize}
   \item One has  
$$ 
\sup_{x\in B(0,R)}\,\sup_{0\le s <t \le T} \frac{\big|V_0(t,x)-V_0(s,x)\big|}{|t-s|^{\kappa_1}} \lesssim (1+R)^\alpha,
$$
   
   \item For all $1\le n \le [p]$ and $1\le k_1,\cdots,k_n \le [p]$, with $\sum_{i=1}^n k_i \le [p]$, for all tuples $I_i \in \{1,\cdots,d\}^{k_i}$, we have
$$
\sup_{x\in B(0,R)}\,\sup_{0\le s <t \le T} \frac{\Big|V_{[I_n]}(t,\cdot)\cdots V_{[I_1]}(t,\cdot)(x)-V_{[I_n]}(s,\cdot)\cdots V_{[I_1]}(s,\cdot)(x)\Big|}{|t-s|^{\kappa_2}} \lesssim (1+R)^\alpha.
$$
\end{itemize}
We assume that the derivative in $x$ of the $V_{[I_n]}(t,\cdot)\cdots V_{[I_1]}(t,\cdot)$ also satisfies the previous estimate.
\end{assum}

\medskip

Let $\bfX$ be an $\bbR^\ell$-valued weak geometric H\"older $p$-rough path. Set ${\bf \Lambda}_{ts} := \log {\bf X}_{ts}$, for all $0\leq s\leq t\leq T$, and denote by $\mu_{ts}$ the time $1$ map of the ordinary differential equation
\begin{equation} 
\label{EqApproximate}
y'_r = (t-s)V_{0}\big(s,y_r(x)\big) + \sum_{k=1}^{[p]} \sum_{I\in\left\{1,\cdots,\ell\right\}^k} \Lambda^{k,I}_{t,s} V_{[I]}\big(s,y_r(x)\big)
\end{equation}
that associates to $x$ the value at time $1$ of the solution path to that equation with initial condition $x$. Note that \textbf{\textsf{Assumption 1}} ensures that equation \eqref{EqApproximate} is well-defined up to time 1 . Following \cite{BailleulRMI}, we define a solution flow to the rough differential equation 
\begin{equation}
\label{EqRDEFlow}
d\varphi_t = V_0(t,\varphi_t)dt + \textrm{F}(t,\varphi_t)d{\bfX}_t,
\end{equation}
where $\textrm{F} := (V_1,\dots,V_\ell)$, as a flow locally well-approximated by $\mu$. Here, we take advantage in this definition of some variant of the definition of \cite{BailleulRMI} introduced by Cass and Weidner in \cite{CassWeidner}. For a parameter $a$, the notation $C_a$ stands for a constant depending only on $a$.

\medskip

\begin{defn}
\label{DefnSolutionFlow}
A flow $\varphi : \Delta_T \times \RR^d \mapsto \RR^d$ is said to be a \textbf{\textsf{solution flow to the rough differential equation}} \eqref{EqRDEFlow} if there exists an exponent $\eta>1$ independent of $\bX$, such that one can associate to any positive radius $R$ two positive constants $C_{R,\bX}$ and $\eps_{\bX}$ such that one has 
\begin{equation}
\label{EqDefnSolFlow}
\sup_{ x\in B(0,R)} \, \big|\varphi_{ts}(x) - \mu_{ts}(x) \big| \le C_{R,\bX}\, |t-s|^{\eta},
\end{equation}
whenever $|t-s|\leq \eps_{\bX}$.
\end{defn}

\medskip

Note that we require the flow to be\textit{ globally defined in time and space}, unlike local flows of possibly exploding ordinary, or rough, differential equations. The latter are only defined on an open set of $\RR_+\times \RR^d$ depending on $\bfX$. This definition differs from the corresponding definition in \cite{BailleulRMI} in the fact that $\epsilon_R$ is required to be independent of $\bfX$. We first state a local in time existence result for the flow, in the spirit of \cite{BailleulRMI}. 

\medskip

\begin{thm}\label{thm:local}
Let the vector fields $V_0$ and $(V_1,\dots,V_\ell)$ satisfy \textbf{\textsf{Assumption 1}} and \textbf{\textsf{Assumption 2}}.   \vspace{0.15cm}

\begin{itemize}
   \item There exists a positive constant $a_1$ such that for all $R>0$, and all $(t,s)\in \Delta_T$ with 
\begin{equation}
\label{EqConditionLocalExistence}
|t-s|^{\frac{1}{p}} (1+R)^{\frac{\alpha}{[p]+1}}  \big(1+\|\bX\|\big) < a_1,
\end{equation}
there is a unique flow $\varphi: [s,t]^2\times B(0,R) \to \bbR^d$ satisfying the estimate \eqref{EqDefnSolFlow} with
\[
\eta=\frac{[p]+1}{p}, \text{ and } C_{R,\bX} = a_2(1+R)^\alpha\big(1+\|\bX\|\big)^{[p]+1},
\]
for some universal positive constant $a_2$. One writes $\varphi(\bfX)$ to emphasize the dependence of $\varphi$ on $\bfX$.   \vspace{0.1cm}

   \item Given a weak geometric rough path $\bfX$ and $(s,t)\in\Delta_T$ and $R$ such that condition \eqref{EqConditionLocalExistence} holds, then $\varphi({\bfX}')$ is well-defined on $[s,t]\times B(0,R)$ for ${\bfX}'$ sufficiently close to $\bfX$, and $\varphi({\bfX}')$ converges to $\varphi(\bfX)$ in $L^\infty\big([s,t]\times B(0,R)\big)$ as ${\bfX}'$ tends to $\bfX$.
\end{itemize}
\end{thm}

\medskip

One says that $\varphi$ depends continuously on $\bfX$ in the topology of uniform convergence on bounded sets. As you can see from the statement of Theorem \ref{thm:local}, the quantity $|t-s|$ is only required in that case to be smaller than a constant depending on $\bfX$ and $R$, unlike what is required from a solution defined globally in time. The proof of Theorem \ref{thm:local} mimics the proof of the analogue local in time result proved in \cite{BailleulRMI}. As the proof of latter contains typos that makes reading it hard, we give in Section \ref{SectionProof} a self-contained proof of this result. 

\medskip

\begin{thm}\label{thm:main}
Let $V_0$ and $(V_1,\dots,V_\ell)$ satisfy \textbf{\textsf{Assumption 1}} and \textbf{\textsf{Assumption 2}}. There exists a unique global in time solution flow $\varphi$ to the rough differential equation \eqref{EqRDEFlow}.   \vspace{0.1cm}
\begin{itemize}
   \item One can choose in the defining relation \eqref{EqDefnSolFlow} for a solution flow 
   $$
   \eta=\frac{1+[p]}{p}, \quad \eps_{\bfX} = c_1\big(1+\|\bX\|\big)^{-p},\quad C_{R,\bfX} = c_2 (1+R)^\alpha \big(1+\|\bX\| \big)^{[p]+1},
   $$ 
   for some universal positive constants $c_1, c_2$.   \vspace{0.1cm}
   
   \item One has for all $f\in C^{[p]+1}_b$ and all $|t-s| \le \eps_{\bfX}$ the estimate
\begin{multline*}
\sup_{x\in B(0,R)} \bigg|f\circ\varphi_{t,s}(x) - \Big\{
 f(x) + (t-s)V_{0}(s,\cdot)f + \sum_{k=1}^{[p]} \sum_{I \in \left\{0,\cdots, \ell\right\}^{k}} X^{k,I}_{t,s}V_{I}(s,\cdot)f\Big\}(x) \bigg|
 \\ \lesssim \|f\|_{C^{[p]+1}_b} (1+R)^{\alpha([p]+1)}\big(1+\|\bfX\|\big)^{[p]+1}|t-s|^{\frac{[p]+1}{p}}.
\end{multline*}   
When $f = \textrm{\emph{Id}}$, one can replace $(1+R)^{\alpha([p]+1)}$ by $(1+R)^\alpha$ and $\|f\|_{C^{n}_b}$ by $1$ in the previous bound.   \vspace{0.1cm}

   \item The map that associates $\varphi$ to $\bfX$ is continuous from the set of weak geometric H\"older $p$-rough paths into the set of continuous flows endowed with the topology of uniform convergence on bounded sets.   \vspace{0.1cm}
   
   \item Finally, there exists two positive universal constants $c_3,c_4$ such that setting 
   $$
   N := \Big[c_3 \big(1+\|\bfX\|\big)^p\Big],
   $$ 
   one has for all $(t,s)\in \Delta_T$,
   $$
   \underset{x\in B(0,R)}{\sup}\, \big|\varphi_{s,t}(x)-x\big| \lesssim 
   \left\{\begin{aligned} 
	   &(1+R)\left( \left(1+c_4 \frac{|t-s|^{\frac1p}N}{(1+R)^{1-\alpha}}\right)^{\frac1{1-\alpha}}-1\right), & \text{if } \alpha<1\\
	   &(1+R)|t-s|^{\frac{1}{p}}e^{c_4 N |t-s|^{\frac1p}}, & \text{if } \alpha=1.
\end{aligned}
\right.
$$
\end{itemize}
\end{thm}

\medskip

The non-trivial part of the proof consists in proving that one can patch together the local flows contructed in Theorem \ref{thm:local} and define a globally well-defined flow. As this requires a careful track of a number of quantities, we provide a proof of the technical results in Section \ref{SectionProof}. Since it is the main contribution of this work, we also give a proof of this theorem using some results of lemmas and propositions of Section \ref{SectionProof}.
 
\medskip 
 
\begin{Dem}[\hspace{-0.03cm}of Theorem \ref{thm:main}] Fix $(s,t)\in\Delta_T$. For $n\ge 0$ and $0\leq k\leq 2^n$ set $t^n_k := k2^{-n}(t-s) + s$ and $\mu^n_{t,s} := \mu_{t^n_{2^n},t^n_{2^n-1}}\circ \cdots \circ \mu_{t^n_{1},t^n_{0}}$. Here is the makor input for the proof of the statement. Proposition \ref{prop:sewing1} below states the  existence of universal positive constants $c_1<1$ and $c_2$ such that for
 \[
 |t-s|^{\frac{1}{p}}\big(1+\|\bfX\|\big) \le c_1^{p} 
 \]
we have for all $n\geq 0$ the estimate
\begin{equation}
\label{eq:mu_dyadic}
\sup_{x\in B(0,R)} \, \big|\mu^n_{ts}(x)-\mu_{ts}(x) \big| \leq c_2 \,|t-s|^{\frac{1+[p]}{p}} \, \big(1+\|\bfX\|\big)^{[p]+1} \, (1+R)^\alpha.
\end{equation}
An elementary Gronwall type bound proved in Lemma \ref{lemma:gronwall_bounds} also gives the estimate
\[
|\mu_{t,s}| \le R + c_2 (1+R)^\alpha.
\]
Putting those two bounds together, one gets the existence of a positive constant $c$ such that one has 
$$
\Big\|\mu_{t^n_k t^n_{k-1}}\circ \cdots \circ \mu_{t^n_1 t^n_{0}}\Big\|_{L^\infty(B(0,R))} \leq R + c\,(1+R)^\alpha,
$$
for all $0\leq k\leq 2^n-1$. Let $n=n(R)$ be the least integer such that
$$
 2^{-n\,\frac{1}{p}} \, (1+R)^{\frac{\alpha}{[p]+1}} \, |t-s|^{\frac{1}p}\big(1+\|\bfX\|\big) \leq \frac{a_1}{(1+2c)^{\frac{\alpha}{[p]+1}}}.
$$
This is the smallest integer such that for all the intervals $(t^n_k,t^n_{k+1})$ satisfy the assumption of Theorem \ref{thm:local}, with starting point $\mu_{t^{n}_{k} t^{n}_{k-1}}\circ \cdots \circ \mu_{t^{n}_{1} t^{n}_{0}}(x)$ and $x\in B(0,R)$. Then, we have for all $m_0,\cdots, m_{2^n-1} \in \bbN$, 
\[
\Big\|\mu^{m_{2^n-1}}_{{t^n_{2^n} t^n_{2^n-1}}}\circ \cdots \circ \mu^{m_{0}}_{t^n_{1} t^n_0} - \mu_{t s} \Big\|_{L^\infty(B(0,R))} \leq c_1\big|t-s\big|^{\frac{1+[p]}{p}}(1+R)^\alpha\big(1+\|\bfX\|\big)^{[p]+1}. 
\]
Sending successively $m_{2^n-1},\dots, m_0$ to $\infty$ and using the continuity of $\varphi$ with respect to its $\RR^d$-valued argument gives 
\begin{equation}\label{eq:homeo}
\Big\|\varphi_{{t^n_{2^n} t^n_{2^n-1}}}\circ \cdots \circ \varphi_{t^n_{1} t^n_0} - \mu_{ts} \Big\| \leq c_2\, \big| t-s \big|^{\frac{1+[p]}{p}} \, (1+R)^\alpha \, \big(1+\|\bfX\|\big)^{[p]+1}. 
\end{equation}
Set, for $x\in B(0,R)$,
$$
\varphi_{ts}(x) := \varphi_{{t^n_{2^n} t^n_{2^n-1}}}\circ \cdots \circ \varphi_{t^n_{1} t^n_0}(x).
$$
Splitting the intervals $(t^n_k,t^n_{k+1})$ into dyadic sub-intervals, one shows that for all $u\in [s,t]$ of the form $u=k2^{-N}(t-s)+s$, one has
\[
\varphi_{t,u}\circ\varphi_{u,s}(x)=\varphi_{t,s}(x).
\]
Finally, since the map
\[
(x,s,t)\to \varphi_{{t^n_{k+1},t^n_{k}}}(x)
\]
is a continuous for all $0\le k\le 2^n-1$, so is $\varphi$. This proves the first item of Theorem \ref{thm:main}.

\medskip

The second item is a byproduct of the bound of Equation \eqref{EqDefnSolFlow} and Corollary \ref{corollary:taylor} below. The third item of the statement is straightforward given that $\varphi$ is constructed from patching together local solution flows.

\medskip

Choose finally a positive constant $c_3$ big enough such that setting 
$$
N := \bigg[c_3\big(1+\|\bfX\|\big)^{p}\bigg],
$$
one has $\frac{t-s}{N} \le \eps_{\bfX}$ and $\big(1+\|\bfX\|\big) N^{-\frac1p} \le 1$. Define also
$$
t_i := \frac{i}{N}(t-s)+s,
$$
and $R_0 := 0$ and 
$$
R_i := \sup_{x\in B(0,R)} \big|\varphi_{t_i s}(x) - x\big|,
$$
for $1\leq i\leq N$. Note that 
$$
C_{|t_{i+1}-t_i|,\|\bfX\|} = \sum_{i=1}^{[p]} \left(\frac{1+\|\bfX\|}{N^{\frac1p}}\right)^{i}|t-s|^{\frac{i}{p} } \lesssim |t-s|^{\frac1p},
$$
for a universal positive multiplicative factor. We thus have 
\begin{equation*}
\begin{split}
\varphi_{t_i s}(x) - x &= \varphi_{t_i t_{i-1}}\big(\varphi_{t_{i-1} s}(x)\big)-\mu_{t_i t_{i-1}}\big(\varphi_{t_{i-1} s}(x)\big)   \\
&\quad+ \mu_{t_i t_{i-1}}\big(\varphi_{t_{i-1} s}(x)\big) - \varphi_{t_{i-1} s}(x)   \\
&\quad+ \varphi_{t_{i-1} s}(x) - x,
\end{split}
\end{equation*} 
and there is an absolute positive constant $K$ such that
$$
R_{i} \le R_{i-1} + K(1+R+R_{i-1})^\alpha|t-s|^{\frac{1}{p}};
$$
the bounds on $\varphi_{ts}(x)-x$ given in the statement follows from that relation.
\end{Dem}

\medskip

As a corollary of Theorem \ref{thm:main}, one proves in Theorem \ref{thm:derivatives_variation} the differentiability of the solution flow with respect to some parameters. This theorem will be of crucial importance in the forthcoming work \cite{BailleulCatellierMeanField}; we state it here in a readily usable form.

\medskip

\begin{assum}
\label{hyp:parameter}
Let $A$ be a Banach, parameter space and let $U$ be a bounded open subset of $A$. Let $(V_i)_{0\le i \leq \ell}$ be time and parameter-dependent vector fields on $\RR^d$ with the following regularity properties.
\begin{itemize}
   \item[$\bullet$] There exists some exponents $\kappa_1>\frac{1+[p]-p}{p}$ and $\kappa_2>\frac{[p]}p$, such that we have
for all integers $\beta_1,\beta_2$ with $0\leq \beta_1+\beta_2 \le [p]+1$,
$$
\sup_{0\le s \le t \le T}\, \frac{\Big\|D_a^{\beta_1}D_x^{\beta_2}V_0(t,.,.)-D_a^{\beta_1}D_x^{\beta_2}V_0(s,.,.)\Big\|_{L^\infty(\RR^d\times U)}}{|t-s|^{\kappa_1}} <+\infty,
$$
   
   \item[$\bullet$] For all $1\le i\leq \ell$, and all integers $\beta_1,\beta_2$ with $0\leq \beta_1+\beta_2 \le [p]+2$, we have
$$
\sup_{0\le s \le t \le T}\, \frac{\Big\|D_a^{\beta_1}D_x^{\beta_2}V_i(t,.,.)-D_a^{\beta_1}D_x^{\beta_2}V_i(s,.,.)\Big\|_{L^\infty(\RR^d\times U)}}{|t-s|^{\kappa_2}} <+\infty
$$ 
\end{itemize}
\end{assum}

\medskip

Refer to Definition \ref{def:accumulation} in  Section \ref{SubsectionControls} for the definition of the {\bf local accumulation} $N_\beta$ of $\bfX$.

\medskip

\begin{thm}\label{thm:derivatives_parameter}
Let $\bfX$ be a $\bbR^\ell$ valued weak geometric Hölder $p$-rough path and suppose that $V_0,V_1,\cdots,V_\ell$ satisfy \textbf{\textsf{Assumptions \ref{hyp:parameter}}}. Let $\varphi(a,\cdot)$ stand for all $a\in U$ for the solution flow to the equation
\begin{equation} \label{eq:flow_parameter}
d\varphi(a,\cdot) = V_{0}\big(t,a,\varphi(a,\cdot)\big)dt + \sigma\big(t,a,\varphi(a,\cdot)\big) d\bfX_t.
\end{equation}
Then for all $1 \le s \le t\le T$, the function $(a,x) \mapsto \varphi_{ts}(a,x)$ is differentiable and   \vspace{0.1cm}

\begin{itemize}
   \item for $|t-s|^{\frac{1}{p}}\big(1+\|\bfX\|\big)\lesssim 1$, and $a\in U$,
$$
\sup_{x \in \RR^d} \, \big|D_a\varphi_{ts}(a,x) \big| \lesssim |t-s|^{\frac1p}\big(1+\| \bfX\|\big)^{[p]}
$$
   
   \item there exists positive constants $\beta$ and $c$ such that one has
$$
\sup_{x \in \RR^d} \, \big|D_a\varphi_{ts}(a,x) \big| \lesssim |t-s|^{\frac1p} \big(1+\|\bfX\|\big)^{}e^{c N_{\beta}}
$$
for all $0\le s \le t \le T$.
\end{itemize}
\end{thm}

\bigskip

\section{Complete proof of Theorem \ref{thm:local}}
\label{SectionProof}

The structure of the proof is simple. One first proves $C^2$ estimates on the time $r$ map of the ordinary differential equation \eqref{EqApproximate}, this is the content of Lemma \ref{lemma:gronwall_bounds}. Building on a Taylor formula given in Lemma \ref{lemma:taylor}, and quantified in Lemma \ref{lemma:bound_remainder} and Corollary \ref{corollary:taylor}, one shows in Proposition \ref{prop:sew} that the $\mu$'s defined what could be called a '\textit{local approximate flow}', after \cite{BailleulRMI}. We then follow the construction recipe of a flow from an approximate flow given in \cite{BailleulRMI}, by patching together the local flows. The crucial global in time existence result is obtained as a consequence of a Gr\"onwall type argument, as can be expected from the fact that, in their simplest form, the growth assumptions of Theorem \ref{thm:main} mean that all the vector fields appearing in the approximate dynamics have $\alpha$-growth. Readers familiar with \cite{BailleulRMI} can go directly to Section \ref{SectionCorollaries}.

\medskip

Recall the definition of $y_r$ as the solution of the ordinary differential equation \eqref{EqApproximate} defining $\mu_{ts}$. The first step in the analysis consists in getting some local in space $C^2$ estimate on $y_r(\cdot)-\textrm{Id}$, with $y_r(\cdot)$ seen as a function of the initial condition $x$ in \eqref{EqApproximate}. Set
$$
C_{|t-s|,\|{\bfX}\|} := |t-s|+\sum_{i=1}^{[p]} |t-s|^{\frac{i}{p}}\|\bfX\|^i.
$$

\medskip

\begin{lem}
\label{lemma:gronwall_bounds}
Assume $V_0$ and $(V_1,\dots,V_\ell)$ satisfy the space regularity \textsf{\textbf{Assumption 1}}, and pick $(s,t)\in\Delta_T$ with 
$$
|t-s|^{\frac{1}{p}}\big(1+\|\bfX\|\big) \leq 1.
$$
Then 
\begin{itemize}
   \item $\big|y_r(x) - x\big| \lesssim \big(1+|x|\big)^\alpha C_{|t-s|,\|{\bfX}\|}$,   \vspace{0.1cm}
   
   \item $\big|Dy_r(x)-\textrm{\emph{Id}}\big| \lesssim C_{|t-s|,\|{\bfX}\|}$,   \vspace{0.1cm}
   
   \item $\big|D^2 y_r(x)\big|\lesssim C_{|t-s|,\|{\bfX}\|}$.   \vspace{0.1cm}
\end{itemize}

The maps $y_r(\cdot)$ are thus $C^1_b$, uniformly in $r\in [0,1]$.
\end{lem}

\medskip

\begin{Dem}
Apply repeatedly Gr\"onwall lemma. We only prove the estimate for $y_r(x)-x$ and leave the remaining details to the reader. It suffices to write
\begin{align*}
|y_r(x) -x| 
	\le& 
		(t-s)\big|V_{0}\big(s,x\big) \big|+ \sum_{k=1}^{[p]} \sum_{I\in\left\{1,\cdots,\ell\right\}^k} \big|\Lambda^{k,I}_{t,s}\big| \, \big|V_{[I]}(s,x)\big|\\
		& +
		(t-s)\int_0^r\Big|V_0\big(s,y_u(x)\big) - V_0(s,x) \Big| \dd u\\
		&+ 
		\sum_{k=1}^{[p]} \sum_{I\in\left\{1,\cdots,\ell\right\}^k} \big|\Lambda^{k,I}_{t,s}\big| \int_0^r \Big|V_{[I]}\big(s,y_u(x)\big) - V_{[I]}(s,x)\Big| \dd u\\
	\lesssim 
		&
		C\big(|t-s|,\|\bfX\|\big) \Big(\big(1+|x|\big)^\alpha +
		\int_0^r |y_u(x) - x| \dd u \Big).
\end{align*}
to get the conclusion from Gr\"onwall lemma, using the fact that $C_{|t-s|,\|{\bfX}\|} \lesssim 1$, for $|t-s|^{\frac{1}{p}}\big(1+\|\bfX\|\big) \leq 1$. The derivative equations satisfied by $Dy_r$ and $D^2y_r$ are used to get the estimates of the statement on these quantities, using once again that the condition of the statement imposes to $C$ to be of order $1$.
\end{Dem}

\medskip

\begin{rem}\label{rem:derivatives1}
Would \textbf{\textsf{Assumption 1}} require in addition that the vector fields $V_0(s,\cdot)$ and $V_{[I]}(s,\cdot)$ were $C^{n+2}_b$ with $\alpha$-growth, uniformly in $0\leq s\leq T$, we would then have the estimate   
$$
\sup_{2\le k \le n+2}\,\big|D^{k}y_r(x)\big| \lesssim C_{|t-s|,\|{\bfX}\|},
$$
under the assumption that $|t-s|^{\frac{1}{p}}\big(1+\|\bfX\|\big) \leq 1$.
\end{rem} 

\medskip

The second step of the analysis is an elementary explicit Taylor expansion; see \cite{BailleulRMI} for the model situation. Given $1\le n\le [p]$, set
$$
\Delta_1^n := \Big\{\big(r_n,\cdots,r_1\big)\in [0,1]^n \ : \ r_n \leq r_{n-1} \leq \cdots \leq r_1\Big\}
$$ 
and 
$$
\fI_{n,[p]} := \Big\{(I_1,\cdots,I_n) \in \left\{1,\cdots,d\right\}^{k_{1}}\times\cdots \times\left\{1,\cdots,d\right\}^{k_{n}} \, ; \sum_{m=1}^n k_m \le [p]\Big\};
$$
indices $k_m$ above are non-null.

\medskip 
 
\begin{lem}\label{lemma:taylor} 
Assume $V_0$ and $(V_1,\dots,V_\ell)$ satisfy the space regularity \textsf{\textbf{Assumption 1}}. For any $1\leq n\leq [p]$ and any vector space valued function $f$ on $\RR^d$ of class $C^n$ we have the Taylor formula
\begin{align*}
 &f\big(\mu_{ts}(x)\big) = 
 f(x) + (t-s)\big(V_{0}(s,\cdot)f\big)(x)   \\ 
 &+ \sum_{i=1}^n \frac1{i!} \sum_{\fI_{i,[p]}} \prod_{m=1}^i \Lambda^{k_m,I_m}_{ts} \, \big(V_{[I_i]}(s,\cdot)\cdots V_{[I_1]}(s,\cdot)f\big)(x)
 \\ 
 &+ \sum_{\fI_{n,[p]}} \prod_{m=1}^n \Lambda^{k_m,I_m}_{ts} \, \int_{\Delta_1^n} \Big\{\big(V_{[I_n]}(s,\cdot)\cdots V_{[I_1]}(s,\cdot)f\big)\big(y_{r_n}(x)\big)   \\
& \qquad\qquad\qquad\qquad\qquad\qquad\qquad\qquad - \big(V_{[I_n]}(s,\cdot)\cdots V_{[I_1]}(s,\cdot)f\big)(x)\Big\}\,dr   \\ 
&+ \eps^{n,f}_{ts}(x)
 \end{align*}
where 
\begin{align*}
&\eps^{n,f}_{ts}(x) := (t-s)\int_0^1\Big\{\big(V_0(s,\cdot)f\big)\big(y_r(x)\big) - \big(V_0(s,\cdot)f\big)(x) \Big\}\,dr   \\ 
&+ \sum_{i=1}^{n-1} \frac{1}{i!} \sum_{\fI_{i,[p]}}(t-s) \prod_{m=1}^{i} \Lambda^{k_m,I_m}_{ts}  \int_{\Delta_1^{i+1}}\big(V_{0}(s,\cdot)V_{[I_i]}(s,\cdot)\cdots V_{[I_1]}(s,\cdot)f\big)\big(y_{r_{i+1}}(x)\big)\,dr   \\ 
&+\sum_{i=2}^{n} \sum_{\substack{\fI_{i-1,[p]}\\k_1+\cdots + k_i \ge [p]+1}} \prod_{m=1}^{i} \Lambda^{k_m,I_m}_{ts}
\int_{\Delta_1^i} \big(V_{[I_i]}(s,\cdot)\cdots V_{[I_1]}(s,\cdot)f\big)\big(y_{r_{i}}(x) \big)\,dr.
\end{align*}
\end{lem}

\medskip

\begin{Dem}
The proof is done by induction, and relies on the following fact. For all $u\in[0,1]$ and all $g\in C^1(B;B)$, we have
$$
g(y_r) - g(x) = (t-s)\int_0^r\big(V_0(s,\cdot)g\big) (y_u)\,du + \sum_{\substack{1\le k\le [p]\\ I \in \left\{1,\cdots,\ell\right\}^k}} \Lambda^{k,I}_{t,s} \int_0^r \big(V_{[I]}(s,\cdot)g\big)(y_u)\,du;
$$
this is step 1 of the induction. For step 2, apply step 1 successively to $g=f$ and $u=1$, then $g=(V_{[I]}(s,\cdot)f\big)$ and $u=r$. This gives
\begin{align*}
f\big(\mu_{ts}(x)\big) &- f(x) = (t-s)\big(V_0(s,\cdot)f\big)(x)   \\
&\quad+ (t-s)\int_0^1 \Big\{ \big(V_0(s,\cdot)f\big)(y_r) -\big(V_0(s,\cdot)f\big)(x)\Big\}\,dr   \\
&\quad+ \sum_{\substack{1\le k\le [p]\\ I \in \left\{1,\cdots,\ell\right\}^k}} \Lambda^{k,I}_{ts} \big(V_{[I]}(s,\cdot)f\big)(x)   \\
&\quad+ \sum_{\substack{1\le k\le [p]\\ I \in \left\{1,\cdots,\ell\right\}^k}} (t-s)\Lambda^{k,I}_{ts} \int_0^1\int_0^{r_1} \big(V_0(s,\cdot)V_{[I]}(s,\cdot)f\big)(y_{r_2})\,dr_2dr_1   \\
&\quad+\sum_{\substack{1\le k_1,k_2\le [p]\\ I_1 \in \left\{1,\cdots,\ell\right\}^{k_1}\\I_2 \in \left\{1,\cdots,\ell\right\}^{k_2}}} \prod_{m=1}^2\Lambda^{k_m,I_m}_{ts} \int_0^1\int_0^{r_1} \big(V_{[I_2]}(s,\cdot)V_{[I_1]}(s,\cdot)f\big)(y_{r_2})\,\,dr_2dr_1.
\end{align*}
The last term of the right hand side can be decomposed into
{\small   \begin{align*}
&\frac12\sum_{\fI_{2,[p]}}\prod_{m=1}^2\Lambda^{k_m,I_m}_{ts}\big(V_{[I_2]}(s,\cdot)V_{[I_1]}(s,\cdot)f\big)(x)   \\
&+ \sum_{\fI_{2,[p]}} \prod_{m=1}^2\Lambda^{k_m,I_m}_{ts} \int_{\Delta_1^2} \Big\{\big(V_{[I_2]}(s,\cdot)V_{[I_1]}(s,\cdot)f\big)(y_{r_2}) - \big(V_{[I_2]}(s,\cdot)V_{[I_1]}(s,\cdot)f\big)(x) \Big\}\,dr_2dr_1   \\
&+ \sum_{\substack{k_1 + k_2 \ge [p] + 1\\ \fI_{1,[p]}}} \prod_{m=1}^2\Lambda^{k_m,I_m}_{ts} \int_{\Delta_1^2} \big(V_{[I_2]}(s,\cdot)V_{[I_1]}(s,\cdot)f\big)(y_{r_2})\,dr_2dr_1;
\end{align*}}
this proves step 2 of the induction. The $n$ to $(n+1)$ induction step is done similarly, and left to the reader.
\end{Dem}

\medskip

Given $f\in C^n_b\big(\RR^d,\RR^d\big)$, set 
$$
\|f\|_n := \big|f(0)\big| + \sup_{k\in\{1,\cdots,n\}}\big\|D^kf\big\|_{\infty}.
$$
A function $g\in C^n\big(\RR^d,\RR^d\big)$ is said to satisfy \textbf{\textsf{Assumption H}} if \textit{for all $1\leq k_1,\dots,k_n\leq [p]$ with $\sum_{i=1}^p k_i\leq [p]$, and all tuples $I_{k_i}\in\{1,\dots,\ell\}^{k_i}$, the functions
$$
V_0(s,\cdot)V_{[I_{n-1}]}(s,\cdot)\cdot V_{[I_1]}(s,\cdot)g\qquad \textrm{and}   \qquad   V_{[I_n]}(s,\cdot)\cdot V_{[I_1]}(s,\cdot)g
$$
are $C^2_b$, with $\alpha$-growth, uniformly in $s\in [0,T]$.   }

\medskip

\begin{lem}\label{lemma:bound_remainder}
Assume \textbf{\textsf{Assumption 1}} holds, and pick a function $f\in C^n_b\big(\RR^d,\RR^d\big)$, for some $2\leq n\leq [p]$. Given $(s,t)\in\Delta_T$ with $|t-s|^{\frac1p}\big(1+\|\bfX\|\big)\leq 1$, we have 
$$
\sup_{x\in B(0,R)}\, \big|\eps^{n,f}_{ts}(x)\big| \lesssim \|f\|_n (1+R)^{n \alpha}\big(1+\|\bfX\|\big)^{[p]+1}(t-s)^{\frac{1+[p]}p},
$$
for all positive radius $R$. 
\begin{itemize}
   \item If furthermore $D^{n+1}f$ exists and is a bounded function, then 
   $$
   \sup_{x\in B(0,R)} \, \Big|D_x\eps^{n,f}_{ts}\Big| \lesssim \|f\|_{n+1}(1+R)^{n\alpha}\big(1+\|\bfX\|\big)^{[p]+1}(t-s)^{\frac{1+[p]}p}.
   $$
   
   \item If finally $f$ satisfies \textsf{\textbf{Assumption H}}, then the previous bound on $D_x\eps^{n,f}_{ts}$ holds with $(1+R)^\alpha$ in place of $(1+R)^{n\alpha}$.
\end{itemize}
\end{lem}

\medskip

\begin{Dem}
Write $dr$ for $dr_i\dots dr_1$ on $\Delta_1^i$, and recall that
\begin{equation}
\label{EqDefnEpsilon}
\begin{split}
&\eps^{n,f}_{ts}(x) = (t-s)\int_0^1 \Big\{\big(V_0(s,\cdot)f\big)(y_r) - \big(V_0(s,\cdot)f\big)(x) \Big\}\,dr   \\ 
&+\sum_{i=1}^{n-1} \frac{1}{i!} \sum_{\fI_{i,[p]}}  (t-s)\prod_{m=1}^{i} \Lambda^{k_m,I_m}_{ts} \int_{\Delta_1^{i+1}} \Big\{V_{0}(s,\cdot)V_{[I_i]}(s,\cdot)\cdots V_{[I_1]}(s,\cdot)f\Big\}(y_{r_{i+1}})\,dr   \\ 
&+\sum_{i=2}^{n} \sum_{\substack{\fI_{i-1,[p]}\\k_1+\cdots + k_i \ge [p]+1}} \prod_{m=1}^{i} \Lambda^{k_m,I_m}_{ts} \int_{\delta_1^i} \Big\{V_{[I_i]}(s,\cdot)\cdots V_{[I_1]}(s,\cdot)f\Big\}(y_{r_{i}})\,dr.
\end{split}
\end{equation}
Recall also that $C_{|t-s|,\|\bfX\|} \lesssim 1$ under the assumption of the statement. 

\ssk

\textcolor{gray}{$\bullet$} As we have for all positive radius $R$, and all points $x,y \in B(0,R)$, the estimate
$$
\big|V_0(s,\cdot)f(x) - V_0(s,\cdot)f(y)\big|\le \|f\|_n \big(1+R\big)^\alpha|x-y|
$$
uniformly in $0\leq s\leq T$, it follows from Lemma \ref{lemma:gronwall_bounds} that
$$
\bigg|(t-s)\int_0^1 \Big\{\big(V_0(s,\cdot)f\big)(y_r) - \big(V_0(s,\cdot)f\big)(x) \Big\}dr\bigg| \lesssim \|f\|_n (t-s) \, C_{|t-s|,\|\bfX\|} \, (1+R)^{2\alpha}.
$$
Note that if $V_0(s,\cdot)f$ is globally Lipschitz continuous one can replace $(1+R)^{2\alpha}$ above by $(1+R)^\alpha$. 

\ssk

We estimate the size of the spatial derivative of the first term in the above decomposition of $\eps^{n,f}_{ts}$ writing
\begin{align*}
\bigg|(t-s)\int_0^1&\Big(D\big(V_0(s,\cdot)f\big)\big(y_r(x)\big) Dy_r(x) - D\big(V_0(s,\cdot)f\big)\big(x\big) \Big)\bigg|\,dr   \\
&\lesssim  
	\bigg|
		(t-s)
		\int_0^1 \dd r 
			\Big(
				D
					\big(V_0(s,\cdot)f\big)
						\big(x\big)
					\Big)
			\Big( 
				Dy_r(x) - \textrm{Id}
			\Big)
	\bigg|   \\
&\quad+
\bigg|
		(t-s)
		\int_0^1 \dd r 
				\Big( D\big(V_0(s,\cdot)f\big)
						\big(y_r(x)\big)
				- D\big(V_0(s,\cdot)f\big)
						\big(x\big)\Big)Dy_r(x)
\bigg|   \\
&\lesssim \|f\|_n\,|t-s|^{1+\frac1p}\, (1+R)^{2\alpha}\,\big(1+\|\bfX\|\big)^{[p]}.
\end{align*}

Once again, one can replace $(1+R)^{2\alpha}$ by $(1+R)^{\alpha}$ if $f$ satisfies \textbf{\textsf{Assumption H}}.

\ssk

\textcolor{gray}{$\bullet$} The two other terms in the decomposition \eqref{EqDefnEpsilon} of $\varepsilon^{n,f}_{ts}$ are estimated in the same way. Remark that 
$$
\sup_{x\in B(0,R)}\Big|\big(V_{0}(s,\cdot)V_{[I_i]}(s,\cdot)\cdots V_{[I_1]}(s,\cdot)f\big)(x)\Big| \lesssim \|f\|_n (1+R)^{(i+1)\alpha}
$$
and that
\begin{align*}
\sup_{x\in B(0,R)}
	\Big| D\big(V_{0}(s,\cdot)V_{[I_i]}(s,\cdot)\cdots V_{[I_1]}(s,\cdot)f\big)(x)\Big|
	\lesssim \|f\|_{n+1} (1+R)^{(i+1)\alpha}.
\end{align*}
One can replace in the previous bounds the first term $(1+R)^{(i+1)\alpha}$ by $(1+R)^\alpha$ and the second term $(1+R)^{(i+1)\alpha}$ by $1$ if $f$ satisfies \textbf{\textsf{Assumption H}}. So
\begin{align*}
\left|\sum_{i=1}^{n-1} \frac{1}{i!} \sum_{\fI_{i,[p]}} (t-s)\right. &\prod_{m=1}^{i} \left.\Lambda^{k_m,I_m}_{ts}  \int_{\Delta_1^{i+1}} \big(V_{0}(s,\cdot)V_{[I_i]}(s,\cdot)\cdots V_{[I_1]}(s,\cdot)f\big)(y_{r_{i+1}}) \,dr \right|   \\
&\lesssim \|f\|_n\, (t-s)^{1+\frac1p}\,(1+R)^{n\alpha}\,\big(1+\|\bfX\|\big)^{[p]}
\end{align*}
and 
\begin{multline*}
\left|\sum_{i=1}^{n-1} \frac{1}{i!} \sum_{\fI_{i,[p]}}   (t-s) \prod_{m=1}^{i} \Lambda^{k_m,I_m}_{ts} \times \right.   \\
\left. \int_{\Delta_1^{i+1}} D\Big\{V_{0}(s,\cdot)V_{[I_i]}(s,\cdot)\cdots V_{[I_1]}(s,\cdot)f\big)\big(y_{r_{i+1}}(x)\Big\}Dy_{r_{i+1}}(x)\,dr\right|   \\
\lesssim \|f\|_{n+1}\, (t-s)^{1+\frac1p}\,(1+R)^{n\alpha}\,\big(1+\|\bfX\|\big)^{[p]}.
\end{multline*}
Once again, if the function $f$ satisfies \textbf{\textsf{Assumption H}}, one can replace $(1+R)^{n\alpha}$ by $(1+R)^\alpha $ in the first bound and $(1+R)^{n\alpha}$ by $1$ in the second bound.

\ssk

The analysis of the last term in the right hand side of the decomposition \eqref{EqDefnEpsilon} for $\varepsilon^{n,f}_{ts}$ is a bit trickier since greater powers of $\|\bfX\|$ can pop out. Indeed, one has
\begin{equation*}
\begin{split}
\bigg|\sum_{i=2}^{n} \sum_{\substack{\fI_{i-1,[p]}\\k_1+\cdots + k_i \ge [p]+1}} &\prod_{m=1}^{i} \Lambda^{k_m,I_m}_{ts} \int_{\Delta_1^i} \Big\{V_{[I_i]}(s,\cdot)\cdots V_{[I_1]}(s,\cdot)f\Big\}(y_{r_{i}})\,dr \bigg|   \\
&\lesssim \|f\|_n\,\sum_{l=1}^{[p]} \big(1+\|\bfX\|\big)^{[p]+i} \, |t-s|^{\frac{[p]+i}{p}} \, (1+R)^{n\alpha}.
\end{split}
\end{equation*}
But recall that $|t-s|^{\frac1p}\big(1+\|\bfX\|\big) \leq 1$, so we have $\big(1+\|\bfX\|\big)^{i-1}|t-s|^{\frac{i-1}{p}} \lesssim 1$,  Hence for all $i\in\{1,\cdots,[p]\}$; this gives the expected upper bound. The same idea is used for the spatial derivatives. Once again, one can replace $(1+R)^{n\alpha}$ by $(1+R)^\alpha$ if the function $f$ satisfies \textbf{\textsf{Assumption H}}.
\end{Dem}

\medskip

\begin{cor}\label{corollary:taylor}
We have
\begin{equation*}
\begin{split}
\sup_{x\in B(0,R)} \bigg|f \circ\mu_{ts}(x) &-  \Big\{ f(x) + (t-s)V_{0}(s,\cdot)f(x) + \sum_{k=1}^{[p]} \sum_{I \in \left\{0,\cdots, \ell\right\}^{k}} X^{k,I}_{ts}V_{I}(s,\cdot)f\Big\}(x) \bigg|   \\ 
&\lesssim \|f\|_{[p]+1}\, (1+R)^{\alpha([p]+1)}\, \big(1+\|\bfX\|\big)^{[p]+1} \, |t-s|^{\frac{[p]+1}{p}}.
\end{split}
\end{equation*}
for all $f\in C^{[p]+1}_b$ with  $\alpha$-growth, and $1 \le k \le [p]$. We also have
\begin{equation*}
\begin{split}
\sup_{x\in B(0,R)} \bigg|\mu_{ts}(x) - \Big(
 x + (t-s)V_{0}(s,x)
 + \sum_{k=1}^{[p]} &\sum_{I \in \left\{0,\cdots, \ell\right\}^{k}} X^{k,I}_{ts}V_{I}(s,x) \Big)\bigg|   \\ 
 &\lesssim  (1+R)^\alpha\big(1+\|\bfX\|\big)^{[p]+1}|t-s|^{\frac{[p]+1}{p}}.
\end{split}
\end{equation*}
and
\begin{equation*}
\begin{split}
\sup_{x\in B(0,R)} \bigg|D\mu_{ts}(x) - \Big( \textrm{\emph{Id}} + (t-s)DV_{0}(s,x) &+ \sum_{k=1}^{[p]} \sum_{I \in \left\{0,\cdots, \ell\right\}^{k}} X^{k,I}_{ts}DV_{I}(s,x) \Big)\bigg|   \\ 
 &\lesssim  (1+R)^\alpha\big(1+\|\bfX\|\big)^{[p]+1}|t-s|^{\frac{[p]+1}{p}}.
\end{split}
\end{equation*}
\end{cor}

\medskip

\begin{Dem}
We only have to bound the sum over $\fI_{n,[p]}$ of the terms
$$
\prod_{m=1}^n \Lambda^{k_m,I_m}_{ts}\int_{\Delta_n}\dd r \Big(\big(V_{[I_n]}(s,\cdot)\cdots V_{[I_1]}(s,\cdot)f\big)(y_{r_n})-\big(V_{[I_n]}(s,\cdot)\cdots V_{[I_1]}(s,\cdot)f\big)(x)\Big)
$$
for $n=[p]$, thanks to lemmas \ref{lemma:taylor} and \ref{lemma:bound_remainder}. We have $k_1 = \cdots = k_{[p]} = 1$, on $\fI_{[p],[p]}$. As we know that we have 
$$
\Big|\big(V_I(s,\cdot)f\big)(x) - \big(V_I(s,\cdot)f\big)(y)\Big| \lesssim (1+R)^{\alpha[p]}\|f\|_{[p]+1}|x-y|,
$$
for all $I\in\left\{1,\cdots,d\right\}^{[p]}$, and all $x,y \in B(0,R)$, it follows that
\begin{equation*}
\begin{split}
\Bigg|\sum_{I_1,\cdots,I_{[p]}\in\left\{1,\cdots,\ell\right\}}
\Big( \Lambda^{1}_{ts}\Big)^{[p]} &\int_{\Delta_1^n} \Big\{\big(V_{I_{[p]}}(s,\cdot)\cdots V_{I_{1}}(s,\cdot) f\big)(y_{r_n})   \\
&\hspace{3cm} -\big(V_{I_{[p]}}(s,\cdot)\cdots,V_{I_{1}}(s,\cdot)f\big)(x)\Big\}\,dr\Bigg|   \\
&\lesssim (1+R)^{\alpha([p]+1)}\,\sum_{k=1}^{[p]}|t-s|^{\frac{i+[p]}{p}} \, \|\bfX\|^{[p]+i}   \\
&\lesssim (1+R)^{\alpha([p]+1)} \, (t-s)^{\frac{1+[p]}{p}} \, \big(1+\|\bfX\|\big)^{[p]+1}.
\end{split}
\end{equation*}
The first estimate of the corollary follows then from the fact that $\exp(\mathbf{\Lambda}_{ts}) = \bfX_{ts}$. The two other estimates are consequences of the fact that the identity map satisfies \textbf{\textsf{Assumption 1}} and \textbf{\textsf{Assumption H}}.
\end{Dem}

\medskip

\begin{rem}\label{rem:derivatives2}
As in \textsf{Remark \ref{rem:derivatives1}}, one can require that $V_0$ and $V_{[I]}$ are more regular, and ask
\begin{enumerate}
   \item[$\bullet$] For all $1\le k_1,\cdots,k_n \le [p]$, with $\sum_{i=1}^n k_i \le [p]$, and all $I_{k_i}\in\{1,\cdots,\ell\}^{k_i}$, the functions
$$
V_0(s,\cdot)V_{[I_{n-1}]}(s,\cdot)\cdots V_{[I_1]}(s,\cdot) \text{ and } V_{[I_{n}]}(s,\cdot)\cdots V_{[I_1]}(s,\cdot)
$$
are $C_b^{2+n}$ with $\alpha$-growth, uniformly in time. 
\end{enumerate}
Under that stronger regularity assumption, we have for all $2\le k \le n+1$,
\begin{equation*}
\begin{split}
\sup_{x\in B(0,R)} \bigg|D^k\mu_{ts}(x) - & \Big\{ (t-s)D^k V_{0}(s,x) + \sum_{j=1}^{[p]} \sum_{I \in \left\{0,\cdots, \ell\right\}^{j}} X^{j,I}_{ts}D^kV_{I}(s,x) \Big\}\bigg|   \\ 
&\lesssim  (1+R)^\alpha \, \big(1+\|\bfX\|\big)^{[p]+1} \, |t-s|^{\frac{[p]+1}{p}}.
\end{split}
\end{equation*}
\end{rem}

\medskip

The next proposition shows that $\mu$ satisfies a localized version of an \textit{approximate flow}; see \cite{BailleulRMI}.

\medskip

\begin{prop}\label{prop:sew}
Given $0\leq s\leq u\leq t\leq T$, with $(t-s)^{\frac{1}{p}}\big(1+\|\bfX\|\big)^{[p]}\leq 1$, we have
\begin{equation*}
\begin{split}
\sup_{x\in B(0,R)}\; &\big|\mu_{tu}\circ \mu_{us}(x)  - \mu_{ts}(x)\big| \vee \bigg|D\big(\mu_{tu}\circ \mu_{us}\big)(x)  - D\big(\mu_{ts}\big)(x)\bigg|   \\
&\lesssim (1+R)^\alpha \,\big(1+\|\bfX\|\big)^{[p]+1} \, |t-s|^{\frac{1+[p]}{p}}.
\end{split}
\end{equation*}
\end{prop}

\ssk

\begin{Dem}
First, remark that
\begin{align*}
\mu_{tu}\circ \mu_{us}(x) &= \mu_{us}(x) + (t-u)V_{0}\big(u,\mu_{us}(x)\big)   \\ 
&\quad+ \sum_{i=1}^{[p]} \frac1{i !} \sum_{\fI_{i,[p]}} \prod_{m=1}^i \Lambda^{k_m,I_m}_{tu} \Big\{V_{[I_i]}(u,\cdot)\cdots V_{[I_1]}(u,\cdot)\Big\}\big(\mu_{us}(x)\big)   \\ 
&\quad+ \tilde \eps_{tu}\big(\mu_{us}(x)\big),
\end{align*}
where $\tilde \eps_{ts}(x) := \eps^{[p],\textrm{Id}}_{ts}(x) + \eps'_{ts}(x)$ and
$$
\eps'_{ts}(x) := \sum_{I\in\{1,\cdots,d\}^{[p]}} \prod_{m=1}^{[p]} \Lambda^{1,i_k}_{t,s} \int_{\Delta_{[p]}} \Big\{(V_I\textrm{Id})(s,y_{r_n}) - (V_I \textrm{Id})(s,x)\Big\}\,dr,
$$
for any $0\leq a\leq b\leq T$. As we also have
\begin{align*}
\mu_{us}(x) - \mu_{ts}(x)
&= -(t-u)V_{0}(s,x)  \\
&\quad+\sum_{i=1}^{[p]} \frac1{i!} \sum_{\fI_{i,[p]}} \Big\{\prod_{m=1}^i \Lambda^{k_m,I_m}_{us}-\prod_{m=1}^i\Lambda^{k_m,I_m}_{ts}\Big\}\big(V_{[I_i]}(s,\cdot)\cdots V_{[I_1]}(s,\cdot)\big)\big(x\big)   \\
&\quad+ \tilde \eps_{us}(x) + \tilde \eps_{ts}(x),
\end{align*}
this gives
{\small  \begin{equation*}
\begin{split}
&\mu_{tu}\circ \mu_{us}(x) - \mu_{ts}(x)   \\  
&= (t-u)\Big(V_{0}\big(u,\mu_{us}(x)\big)-V_{0}\big(s,\mu_{us}(x)\big)\Big)   \\
&\quad+ (t-u)\Big(V_{0}\big(s,\mu_{us}(x)\big)-V_{0}\big(s,x\big)\Big)   \\
&\quad+ \sum_{i=1}^{[p]} \frac1{i!} \sum_{\fI_{i,[p]}} \Big\{\prod_{m=1}^i \Lambda^{k_m,I_m}_{tu}+\prod_{m=1}^i \Lambda^{k_m,I_m}_{us} - \prod_{m=1}^i\Lambda^{k_m,I_m}_{ts}\Big\} \big(V_{[I_i]}(s,\cdot)\cdots V_{[I_1]}(s,\cdot)\textrm{Id}\big)(x)   \\
&\quad+ \sum_{i=1}^{[p]} \frac1{i!} \sum_{\fI_{i,[p]}} \prod_{m=1}^i \Lambda^{k_m,I_m}_{tu}  \Big\{\big(V_{[I_i]}(u,\cdot)\cdots V_{[I_1]}(u,\cdot)\textrm{Id}\big)\big(\mu_{us}(x)\big)   \\ 
&\hspace{6cm} - \big(V_{[I_i]}(s,\cdot)\cdots V_{[I_1]}(s,\cdot)\textrm{Id}\big)\big(\mu_{us}(x)\big) \Big\}   \\
&\quad+ \sum_{i=1}^{[p]} \frac1{i !}\sum_{\fI_{i,[p]}} \prod_{m=1}^i \Lambda^{k_m,I_m}_{tu} \Big\{\big(V_{[I_i]}(s,\cdot)\cdots V_{[I_1]}(s,\cdot)\textrm{Id}\big)\big(\mu_{u,s}(x)\big)   \\
&\hspace{6cm}- \big(V_{[I_i]}(s,\cdot)\cdots V_{[I_1]}(s,\cdot)\textrm{Id}\big)(x) \Big\}   \\
&\quad+ \tilde \eps_{tu}\big(\mu_{us}(x)\big) + \tilde \eps_{us}\big(x\big) + \tilde \eps_{ts}(x)   \\   
&=: (1) + \cdots + (6).
\end{split}
\end{equation*}   }
The bounds of the statement can be read on that decomposition; we give the details for $(\mu_{tu}\circ\mu_{us})(x)$ and live the details of the estimate for its derivative to the reader.

\ssk

It follows from \textbf{\textsf{Assumption 2}} on the time regularity of $V_0$ and the $V_{[I_i]}\dots V_{[I_1]}\textrm{Id}$ that 
\begin{align*}
\big|(1)\big| + \big| (4)\big| &\lesssim (1+R)^\alpha \Big((t-u)(u-s)^{\kappa_1} +\big(1+\|\bfX\|\big)^{[p]}(t-u)^{\frac1p}(u-s)^{\kappa_2}\Big)\\
 &\lesssim (1+R)^\alpha (t-s)^{\frac{1+[p]}{p}} \big(1+\|\bfX\|\big)^{[p]}. 
 \end{align*}
Lemma \ref{lemma:bound_remainder} takes care of the remainder terms $(6)$. By using lemma \ref{lemma:gronwall_bounds} and the fact that $V_0$ is Lipschitz continuous in space, uniformly in time, one gets
$$
\big| (2) \big| \lesssim (t-u)(u-s)^{\frac1p}\big(1+\|\bfX\|\big)^{[p]}(1+R)^\alpha \lesssim (t-s)^{\frac{1+[p]}p}\big(1+\|\bfX\|\big)^{[p]}(1+R)^\alpha.
$$
To estimate the terms $(3)$ and $(5)$, set 
$$
g(s,\cdot) := V_{[I_i]}(s,\cdot)\cdots V_{[I_1]}(s,\cdot)\textrm{Id}.
$$
We start by doing a Taylor expansion of $g\big(s,\mu_{ts}(x)\big)$ using Lemma \ref{lemma:taylor}, to the order $n=[p]-\sum_{j=1}^i k_j$. As $g(s,\cdot)$ satisfies \textbf{\textsf{Assumption H}} as a consequence of \textbf{\textsf{Assumption 1}}, one can use Lemma \ref{lemma:bound_remainder} to get the expected bounds, using the fact that $\bfX_{u,s} \bfX_{t,u} = \bfX_{t,s}$ and $\exp({\bf \Lambda}) = \bfX$. Details of these algebraic computations can be found in the proof of the corresponding statement in \cite{BailleulRMI}.
\end{Dem}

\medskip

\begin{Remark*}
One has similar local bounds for higher derivatives of $\mu_{tu}\circ\mu_{us}-\mu_{ts}$ in the setting of \emph{\textsf{Remark \ref{rem:derivatives2}}}.
\end{Remark*}

\medskip

Write here part of the conclusion of Proposition \ref{prop:sew} under the form 
$$
\sup_{x\in B(0,R)}\; \big|\mu_{t,u}\circ \mu_{u,s}(x)  - \mu_{t,s}(x)\big| \leq C_0\, (1+R)^\alpha \,\big(1+\|\bfX\|\big)^{[p]+1} \, |t-s|^{\frac{1+[p]}{p}},
$$
for some positive constant $C_0$. Given $n\geq 1$, and $0\leq s\leq t\leq T$, set $t^n_k := k2^{-n}(t-s) + s$. Pick $\epsilon_0$ such that 
$$
2^{-\frac{1+[p]-p}{p}}(1+2\epsilon_0) < 1
$$
and 
$$
L > \frac{C_0}{1-2^{-\frac{1+[p]-p}{p}}(1+2\epsilon_0)}.
$$

\medskip

\begin{prop}\label{prop:sewing1}
For all $0\le s<t\le T$ with 
$$
L\,|t-s|^{\frac{1}{p}} \, \big(1+\|\bfX\|\big)\leq \epsilon_0,
$$
and all positive radius $R$, we have
$$
\sup_{|x|\leq R}\;\Big|\mu_{t^{n}_{2^{n}}t^{n}_{2^{n}-1}}\circ \cdots \circ \mu_{t^{n}_1 t^n_0}(x) - \mu_{ts}(x)\Big| \leq L\, |t-s|^{\frac{1+[p]}{p}} \,(1+R)^\alpha \, \big(1+\|\bfX\|\big)^{[p]+1}.
$$
\end{prop}

\medskip

\begin{Dem}
The proof is done by induction on $n$. Note first that we can take $L$ enough to have $\frac{\epsilon_0}{L} \leq 1$ and $C_{|t-s|,\|\bfX\|} \lesssim \epsilon_0$. Proposition \ref{prop:sew} provides the initialisation of the induction. Assume step $n$ of the induction has been proved and set 
$$
u := \frac{t+s}{2} = t^{n+1}_{2^n},
$$ 
so the statement of the proposition holds on the intervals $(s,u)$ and $(u,t)$. We have
\begin{equation*}
\begin{split}
\Big|\big(\mu_{t^{n+1}_{2^{n}}t^{n+1}_{2^{n}-1}}\circ \cdots \circ &\mu_{t^{n+1}_{1} t^{n+1}_{0}}\big)(x)\Big|   \\
&\leq \Big|\big(\mu_{t^{n+1}_{2^{n}} t^{n+1}_{2^{n}-1}}\circ \cdots \circ \mu_{t^{n+1}_{1} t^{n+1}_{0}}\big)(x)- \mu_{us}(x)\Big| + \big|\mu_{us}(x)\big|   \\ 
&\leq   L \, 2^{-\frac{1+[p]}{p}} |t-s|^{\frac{1+[p]}{p}} \, (1+R)^\alpha \, \big(1+\|\bfX\|\big)^{[p]+1} + R   \\
&\quad + C_{|t-s|,\|\bfX\|} \, (1+R)^\alpha   \\
&\leq R + 2(1+R)^\alpha \, \epsilon_0.
\end{split}
\end{equation*}
and 
$$
\sup_{x\in B(0,R)}\,\big| D\mu_{tu}(x)\big| \leq 1+2\epsilon_0,
$$
by Lemma \ref{lemma:gronwall_bounds}. Furthermore we have
\begin{equation*}
\begin{split}
&\mu_{t^{n+1}_{2^{n+1}}t^{n+1}_{2^{n+1}-1}}\circ \cdots \circ \mu_{t^{n+1}_{1}t^{n+1}_{0}}(x) - \mu_{ts}(x)   \\
&= \Big( \mu_{t^{n+1}_{2^{n+1}} t^{n+1}_{2^{n+1}-1}}\circ \cdots \circ \mu_{t^{n+1}_{2^n+1} t^{n+1}_{2^n}} - \mu_{tu} \Big)
\circ\Big(\mu_{t^{n+1}_{2^{n}} t^{n+1}_{2^{n}-1}}\circ \cdots \circ \mu_{t^{n+1}_{1} t^{n+1}_{0}}\Big)(x)   \\
&\quad+ \mu_{tu} \circ\big(\mu_{t^{n+1}_{2^{n}} t^{n+1}_{2^{n}-1}}\circ \cdots \circ \mu_{t^{n+1}_{1} t^{n+1}_{0}}\big)(x) - \mu_{tu}\circ\mu_{us}(x)   \\
&\quad+\mu_{tu}\circ\mu_{us}(x)-\mu_{ts}(x).
\end{split}
\end{equation*}
We thus have for all $x\in B(0,R)$, the estimate
\begin{align*}
\Big|\mu_{t^{n+1}_{2^{n+1}} t^{n+1}_{2^{n+1}-1}}\circ \cdots &\circ \mu_{t^{n+1}_{1} t^{n+1}_{0}}(x) - \mu_{ts}(x)\Big|   \\
&\leq L \, \left|\frac{t-s}{2}\right|^{\frac{1+[p]}{p}}\Big(1+R + 2\eps_0(1+R)^\alpha\Big)^\alpha\big(1+\|\bfX\|\big)^{[p]+1}   \\
&+ (1+2\eps_0)\,L \, (1+R)^\alpha \left|\frac{t-s}{2}\right|^{\frac{1+[p]}{p}}\big(1+\|\bfX\|\big)^{[p]+1}   \\
&+ C_0\,|t-s|^\frac{1+[p]}{p}(1+R)^\alpha \big(1+\|\bfX\|\big)^{[p]+1},
\end{align*}
from which the induction step follows given our choice of $\epsilon_0$ and $L$.
\end{Dem}

\medskip

The same bound for the derivative of the approximate flow requires a bound on $|t-s|$ that \emph{depends on $(1+R)^\alpha$}, such as described here. 

\medskip

\begin{prop}\label{prop:sewing2}
One can find a positive constant $\epsilon_1<1$ such that for $0\leq s\leq t\leq T$ with 
$$
C_0|t-s|^{\frac{1}{p}}(1+R)^{\frac\alpha{1+[p]}}\big(1+\|\bfX\|\big)\leq \epsilon_1,
$$
we have, for all positive radius $R$,
$$
\sup_{|x|\le R}\Big|D\big(\mu_{t^{n}_{2^{n}} t^{n}_{2^{n}-1}}\circ \cdots \circ \mu_{t^{n}_{1} t^{n}_{0}}\big)(x) - D\big(\mu_{ts}\big)(x)\Big| \leq L\,  |t-s|^{\frac{1+[p]}{p}} \, (1+R)^\alpha \, \big( 1+\|\bfX\|\big)^{[p]+1}.
$$
\end{prop}

\medskip

\begin{Dem}
The proof is a variation on the theme of the proof of Proposition \ref{prop:sewing1}. We provide the details for the reader's convenience, and keep the notation $u$ for $\frac{s+t}{2}$. We proceed here as well by induction and loot at the '$n$ to $n+1$' induction step of the proof.

\begin{align*}
D&\Big(\mu_{t^{n+1}_{2^{n+1}} t^{n+1}_{2^{n+1}-1}}\circ \cdots \circ \mu_{t^{n+1}_{1} t^{n+1}_{0}} \Big)(x)- D\mu_{ts}(x)   \\
&= \Big\{D\big(\mu_{t^{n+1}_{2^{n+1}} t^{n+1}_{2^{n+1}-1}}\circ \cdots \circ \mu_{t^{n+1}_{2^n+1} t^{n+1}_{2^n}}\big) - D\mu_{tu} \Big\} \big(\mu_{t^{n+1}_{2^{n}} t^{n+1}_{2^{n}-1}}\circ \cdots \circ \mu_{t^{n+1}_{1} t^{n+1}_{0}}\big)(x)   \\
&\qquad\qquad \times D\big(\mu_{t^{n+1}_{2^{n}}t^{n+1}_{2^{n}-1}}\circ \cdots \circ \mu_{t^{n+1}_{1}t^{n+1}_{0}}\big)(x)   \\
&\quad+ \bigg\{\big(D\mu_{tu}\big) \Big(\big(\mu_{t^{n+1}_{2^{n}} t^{n+1}_{2^{n}-1}}\circ \cdots \circ \mu_{t^{n+1}_{1} t^{n+1}_{0}}\big)(x)\Big)-D\mu_{tu}\big(\mu_{us}(x)\big)\bigg\}   \\
&\qquad\qquad \times D\big(\mu_{t^{n+1}_{2^{n}},t^{n+1}_{2^{n}-1}}\circ \cdots \circ \mu_{t^{n+1}_{1},t^{n+1}_{0}}\big)(x)   \\
&\quad+ D\mu_{tu}\big( \mu_{us}(x)\big)\Big(D\big(\mu_{t^{n+1}_{2^{n}} t^{n+1}_{2^{n}-1}}\circ \cdots \circ \mu_{t^{n+1}_{1} t^{n+1}_{0}}\big)(x)-D\mu_{us}(x)\Big)   \\
&\quad+ D\big(\mu_{tu}\circ\mu_{us}\big)(x) - D\mu_{ts}(x).
\end{align*}
We know from the induction step and the $R$-dependent assumption on $u-s$ that 
\begin{align*}
\Big|D\Big(\mu_{t^{n+1}_{2^{n}} t^{n+1}_{2^{n}-1}}\circ \cdots \circ \mu_{t^{n+1}_{1} t^{n+1}_{0}}\Big)(x) - D\mu_{us}(x)\Big| \leq
\epsilon_1 \, 2^{-\frac{1+[p]}{p}},
\end{align*}
and
$$
\Big|D\Big(\mu_{t^{n+1}_{2^{n}} t^{n+1}_{2^{n}-1}}\circ \cdots \circ \mu_{t^{n+1}_{1} t^{n+1}_{0}}\Big)(x)\Big| \leq 1+2\epsilon_1.
$$
We also have from Proposition \ref{prop:sewing1}
$$
\Big|\Big(\mu_{t^{n+1}_{2^{n}}t^{n+1}_{2^{n}-1}}\circ \cdots \circ \mu_{t^{n+1}_{1} t^{n+1}_{0}}\Big)(x)\big| \leq \big(1+2\epsilon_1\big)\,(1+R)-1,
$$
and Lemma \ref{lemma:gronwall_bounds} gives us a uniform control on the Lipschitz size of the $\mu_{ba}$. We thus have
{\small   \begin{equation*}
\begin{split}
\Big|&D\Big(\mu_{t^{n+1}_{2^{n+1}} t^{n+1}_{2^{n+1}-1}}\circ \cdots \circ \mu_{t^{n+1}_{1} t^{n+1}_{0}} \Big)(x)- D\mu_{ts}(x)\Big|   \\
&\leq  L\,(1+2\epsilon_1)^{1+\alpha} \, 2^{-\frac{1+[p]}{p}} \, |t-s|^{\frac{1+[p]}{p}} \, \big(1+\|\bfX\|\big)^{[p]+1} \, (1+R)^\alpha   \\
&\quad+ L\,2^{-\frac{1+[p]}{p}}\, \epsilon_1\, (1+2\epsilon_1) \, |t-s|^{\frac{1+[p]}{p}} \, \big(1+\|\bfX\|\big)^{[p]+1} \, (1+R)^\alpha   \\
&\quad+ L\,( 1 + \epsilon_1 ) \, 2^{-\frac{1+[p]}{p}} \, |t-s|^{\frac{1+[p]}{p}} \, \big(1+\|\bfX\|\big)^{[p]+1} \, (1+R)^\alpha   \\
&\quad+ C_0\,|t-s|^{\frac{1+[p]}{p}} \,\big(1+\|\bfX\|\big)^{[p]+1} \, (1+R)^\alpha   \\
&\leq \Big(C_0+b 2^{-\frac{1+[p]}{p}}\big((1+2\eps)^\alpha + \epsilon_1(1+2\epsilon_1) + (1+\epsilon_1)\big)\Big)\, |t-s|^{\frac{1+[p]}{p}} \, (1+R)^\alpha \, \big( 1+\|\bfX\|\big)^{[p]+1}.
\end{split}
\end{equation*}   }
An adequate choice of $\epsilon_1$ closes the induction step, given the definition of $L$.
\end{Dem}

\medskip

\begin{rem}\label{rem:derivatives3}
In the improved regularity conditions on the vector fields stated in \emph{Remark \ref{rem:derivatives2}}, we have for all $2 \le k \le n+1$ and for all 
$$
C_0\,|t-s|^{\frac{1}{p}} \, (1+R)^{\frac{k\alpha}{1+[p]}} \, \big(1+\|\bfX\|\big)\leq \eps_1,
$$
one have
$$
\sup_{|x|\le R}\Big|D^k\big(\mu_{t^{n}_{2^{n}} t^{n}_{2^{n}-1}}\circ \cdots \circ \mu_{t^{n}_{1} t^{n}_{0}}\big)(x) - D^k\mu_{ts}(x)\Big| \leq L \, |t-s|^{\frac{1+[p]}{p}} \, (1+R)^{\alpha} \, \big( 1+\|\bfX\|\big)^{[p]+1}.
$$
\end{rem}

\medskip

With all these preliminary results at hand, we are now in a position to give a proof of our local well-posedness result, Theorem \ref{thm:local}.

\medskip

\begin{Dem}[of Theorem \ref{thm:local}]\label{page}
We treat existence and uniqueness one after the other. We keep the above notations, and set, in addition,
$$
\mu^n_{ts}=\mu_{t^n_{2^n} t^n_{2^n-1}}\circ\cdots \circ\mu_{t^n_1 t^n_0}.
$$

\ssk

\noindent \textbf{\textsf{Local in time existence -- }} For all $x \in B(0,R)$, 
\begin{multline}\label{eq:extension_derivatives}
		\mu^{n+1}_{ts}(x)-\mu^{n}_{ts}(x) = \\
			\sum_{k=1}^{2^n}\Big(
				\mu_{t^{n+1}_{2^{n+1}} t^{n+1}_{2^{n+1}-1}}
					\circ \cdots \circ
				\mu_{t^{n+1}_{2k+3} t^{n+1}_{2k+2}} 
					\circ \big(\mu_{t^{n+1}_{2k+2} t^{n+1}_{2k+1}}\circ \mu_{t^{n+1}_{2k+1} t^{n+1}_{2k}}\big)
				\\
					-
					\mu_{t^{n+1}_{2^{n+1}} t^{n+1}_{2^{n+1}-1}}
					\circ \cdots \circ
				\mu_{t^{n+1}_{2k+3} t^{n+1}_{2k+2}} 
					\circ \big(\mu_{t^{n+1}_{k+1} t^{n+1}_{2k}}
				\big)
					\Big)\\
					\circ 
				\mu_{t^{n}_{k} t^{n}_{k-1}} 
					\circ \cdots \circ 
				\mu_{t^{n}_1t^{n}_0}(x).
\end{multline}
It follows from Proposition \ref{prop:sewing2} that the maps 
$$
\mu_{t^{n+1}_{2^{n+1}} t^{n+1}_{2^{n+1}-1}} \circ \cdots \circ \mu_{t^{n+1}_{2k+3} t^{n+1}_{2k+2}}
$$ 
are Lipschitz continuous, uniformly in $n$, with a Lipschitz constant that depends neither on $\bfX$ nor on $R$.  Furthermore, thanks to Proposition \ref{prop:sewing1}, 
$$
\Big|\mu_{t^{n}_{k} t^{n}_{k-1}} \circ \cdots \circ \mu_{t^{n}_{1} t^{n}_{0}}(x)\Big| \le R + 2\epsilon_1 \,(1+R)^\alpha.
$$ 
Finally, Proposition \ref{prop:sew} tells us that
\begin{equation}\label{eq:var1}
\big|\mu^{n+1}_{ts}(x)-\mu^{n}_{ts}(x)\big| \lesssim 2^{-n\frac{[p]+1-p}{p}} \, |t-s|^{\frac{1+[p]}{p}} \, (1+R)^\alpha \, \big(1+\|\bfX\|\big)^{[p]+1}.
\end{equation}
The sequence $\mu^n_{ts}$ is thus uniformly convergent on the ball $B(0,R)$ to a limit, continuous, function denoted by $\varphi_{ts}$; it satisfies the estimate 
$$
\sup_{x\in B(0,R)}\big|\varphi_{ts}(x) - \mu_{ts}(x)\big| \lesssim |t-s|^{\frac{1+[p]}{p}}\, (1+R)^\alpha \, \big(1+\|\bfX\|\big)^{[p]+1}.
$$
Finally, for all dyadic points $a\in[s,t]$ and all $x\in B(0,R)$, we have by construction
\begin{align*}
\varphi_{ta}(x)\circ \varphi_{as}(x) = \varphi_{ts}(x).
\end{align*}
As $\bfX$ is an Hölder continuous rough path, the function $(x; s,t)\mapsto \mu_{ts}(x)$, from $B(0,R)\times \{0 \le s < t \le T\}$ to $\bbR^d$, is continuous. The continuity of $\varphi$ as a function of $(x;s,t)$ follows in a straightforward way; its continuous dependence on $\bfX$ is a consequence of the continuous dependence of $\mu$ with respect to $\bfX$. Note however that $\varphi_{ts}$ is only defined at that stage for $s$ and $t$ close enough.
 
\medskip
 	
\noindent \textsf{\textbf{Uniqueness -- }} Let $\psi$ stand for another solution flow, with associated constants $\overline{\epsilon}_{\bfX}$ and $\overline{C}_{R,\bfX}$, and exponent $\overline{\eta}>1$. Take $R$ and $(s,t)$ satisfying the conditions of Proposition \ref{prop:sewing2}, with $|t-s|\leq \overline{C}_{R,\bfX}$. Then 
\begin{equation}\label{eq:local_uniqueness}
\begin{aligned}
\big|&\varphi_{ts}(x) - \psi_{ts}(x)\big| \leq \Big|\mu^n_{ts}(x)-\varphi_{ts}(x)\Big| + \Big|\mu^n_{ts}(x)- \psi_{ts}(x)\Big|   \\
&\leq \Big|\mu^n_{ts}(x)-\varphi_{ts}(x)\Big|   \\
&\quad+ 
	\sum_{k=0}^{2^n-1}
		\bigg|
			\Big(\big(\mu_{t^n_{2^n} t^n_{2^n-1}}\circ \mu_{t^n_{k+2} t^n_{k+1}} \big)\circ \mu_{t^n_{k+1} t^n_{k}} \\
&\quad \quad \quad \quad \quad			- \big(\mu_{t^n_{2^n} t^n_{2^n-1}}\circ \mu_{t^n_{k+2} t^n_{k+1}}\big) \circ \psi_{t^n_{k+1} t^n_{k}}
			\Big)
				\circ  \psi_{t^n_{k} t^n_{k-1}}\circ \cdots \circ \psi_{t^n_{1} t^n_{0}}(x)
		\bigg|   \\
&\lesssim \big|\mu^n_{ts}(x)-\varphi_{ts}(x)\Big| + 2^{-n(\overline{\eta} - 1)}.
\end{aligned}
\end{equation}
Local uniqueness follows from that estimate. We have used here the fact that the $\mu^n_{ts}$ are Lipschitz continuous, uniformly in $n$, and that
$$
\sup_{x\in B(0,R)}\, \big|\psi_{t,s}-\mu_{t,s}(x)\big| \lesssim \overline{C}_{R,\bfX} \; |t-s|^{\overline{\eta}}.
$$
\end{Dem}

\bigskip
\bigskip

\section{Corollaries and extensions}
\label{SectionCorollaries}

We emphasize in the Section \ref{SubsectionYoung} and Section \ref{SubsectionDerivativeFlow} two consequences on solutions to rough differential equations of the above results/computations. Young and mixed rough/Young equations are considered in Section \ref{SubsectionYoung}, and differentiability of the solution flow with respect to parameters is considered in Section \ref{SubsectionDerivativeFlow}. The estimates on the derivative flow we get there will be used in the forthcoming work \cite{BailleulCatellierMeanField} on limit theorems for systems of mean field rough differential equations. We worked so far in with weak geometric H\"older $p$-rough paths; one can actually work with general rough paths, controlled by arbitrary controls \cite{Lyons98}. A non-explosion criterion with quantitative estimates is provided in Section \ref{SubsectionControls} in this more general setting.

\bigskip

\subsection{Young and mixed rough-Young differential equations}
\label{SubsectionYoung}

The proofs of theorems \ref{thm:local} and \ref{thm:main} do not use the fact the the drift term is driven only by time. Instead we treat the signal $t\to t$ as a Lipschitz path, and deal with it using Young differential calculus techniques. A direct counterpart of this approach is a loss of regularity in the coefficients, either in time and space. A real reward of this approach, which does not modify the proof but requires only more notations, is an extension of the results to a mixed Young-Rough differential equation. 

Let $V_0$ and $\textrm{F} = (V_1,\dots,V_\ell)$ be given; let another family $\textrm{G} := (W_1,\dots, W_m)$ of vector fields on $B$ be given. A solution flow to the mixed rough-Young differential equation is defined as in Definition \ref{DefnSolutionFlow}, with the 'approximate flow'  $\mu_{ts}$ defined as the time $1$ map of the ordinary differential equation
$$
y'_r = V_0\big(s,y_r\big) (t-s) + \sum_{j=1}^m Y^j_{ts}\,W_j\big(s,y_r\big) + \sum_{k=}^{[p]} \sum_{I\in \{1,\cdots,\ell\}^k} \Lambda^{k,I}_{ts}\, V^X_{[I]} \big(s,y_r\big).
$$
The constants $\varepsilon$ and $C$ that appear in the defining estimate \eqref{EqDefnSolFlow} are now allowed to depend on $R,\bfX$ and $Y$.

\medskip

\begin{cor}
Let $\bfX$ be an $\bbR^\ell$-valued weak geometric Hölder $p$-rough path and $Y$ be an $\bbR^m$-valued $\frac{1}{q}$-H\"older path, with $\frac1p+\frac1q >1$ and $p\geq 2$. Assume $(V_0,\textrm{\emph{F}})$ and $(W_i,\textrm{\emph{F}})$  satisfy \textbf{\textsf{Assumption 1}} and \textbf{\textsf{Assumption 2}} for all $1\leq i\leq m$. Assume furthermore that there exists a positive exponent $\kappa$ such that $\kappa+\frac{1}{q}>1$, and
$$
\sup_{x\in B(0,R)} \, \sup_{0\le s < t \le T} \, \frac{\big|W_{i}(t,x) - W_{i}(s,x) \big|}{|t-s|^\kappa} \lesssim (1+R)^\alpha.
$$
Then the rough differential equation
$$
d\varphi_t = V_0(t,\varphi_t) dt + \textrm{\emph{G}}(t,\varphi_t) \,dY_t + \textrm{\emph{F}}(t,\varphi) \,d\bfX_t,
$$
has a unique global in time solution flow.
\end{cor}

\medskip

On can choose the constants $\varepsilon_{{\bfX},Y}$ and $C_{R,{\bfX},Y}$ such that 
$$
(t-s)^\frac{1}{q}\,\big(1+ \|Y\|\big) + (t-s)^\frac{1}{p}\big(1+ \|\bfX\|\big) \lesssim 1,
$$
and 
$$
C_{R,{\bfX},Y} \simeq (1+R)^\alpha\,\big(1+ \|Y\| + \|\bfX\|\big)\,\big(1+ \|\bfX\|\big)^{[p]}
$$
and 
$$
N \simeq \max\Big\{\big[(1+\|Y\|)^{q}\big], \big[(1+\|\bfX\|)^{p}\big]\Big\}.
$$

The proof is left to the reader since it is a direct modification of Section \ref{SectionProof}, with more notations.

\bigskip

\subsection{Derivative flow}
\label{SubsectionDerivativeFlow}

Rough differential equations
\begin{equation}
\label{EqRDEDerivativeFlow}
d\varphi_t = V_0(t,\varphi_t) dt + \textrm{F}(t,\varphi_t)\,d{\bf}_t
\end{equation}
generate flows of diffeomorphisms under appropriate regularity conditions on the driving vector fields. The pair $(\varphi,D\varphi)$, made up of $\varphi$ and its differential, also satisfies an equation, with 'triangular' structure
$$
d(D\varphi) = DV_0(t,\varphi_t)D\varphi_t \, dt + D\textrm{F}(t,\varphi_t)D\varphi_t\,d\bfX_t.
$$
One can find results on derivative flows in the book \cite{FrizVictoirBook} of Friz and Victoir, Chapter 11; see also the interesting works \cite{LejayCoutinOmega} and \cite{CoutinLejayPerturbed} of Coutin and Lejay.  One gets another proof of the differentiability of the flow with respect to the initial point as a direct byproduct of the results of Section \ref{SectionFlows}. Pick $p>2$.

\medskip

\begin{assum}\label{hyp:derivatives}
Let $V_0$ and $V_1,\dots, V_\ell$ be a set of time dependent vector fields on $B$ such that there exists two exponents with $\kappa_1 >\frac{1+[p]-p}{p}$, and $\kappa_2 + \frac{1}{p} > 1$, such that
$$
\sup_{0\le s < t \le T}\,\frac{\big\|V_0(t,\cdot) - V_0(s,\cdot)\big\|_{C^{2+n}_b}}{|t-s|^{\kappa_1}} < + \infty
$$
and each $V_i$ satisfies the estimate
$$
\sup_{0\le s < t \le T}\, \frac{\big\|V_i(t,\cdot) - V_i(s,\cdot)\big\|_{C^{3+n}_b}}{|t-s|^{\kappa_2}} < + \infty.
$$
\end{assum}

\medskip

\begin{thm}\label{thm:derivatives}
Let $\bfX$ be a weak geometric Hölder  $p$-rough path and $(V_0,V_1,\cdots,V_\ell)$ which satisfy \textbf{\textsf{Assumption \ref{hyp:derivatives}}}. Let $\varphi$ stand for the solution flow to the rough differential equation \eqref{EqRDEDerivativeFlow}. Then each $\varphi_{ts}$ is of class $C^n$, has linear growth and bounded derivatives. Furthermore for a suitable positive constant $\epsilon_3$, independent of $\bfX$, and $|t-s|^{\frac{1}{p}}\big(1+\|\bfX\|\big) \leq \epsilon_3$, we have
$$
\Big\| D^k\varphi_{ts} -D^k \mu_{ts}\Big\|_\infty \lesssim |t-s|^{\frac{1+[p]}{p}}\big(1+\|\bfX\|\big)^{[p]+1},
$$
for all $0\leq k\leq n$. Finally there exists some positive constants $c_1,\cdots,c_n$, independent of $\bfX$, such that \emph{for all} $0\le s \le t \le T$, and every $1\leq k\leq n$, we have
$$
\sup_{x\in \bbR^d} \big| \varphi_{ts}(x) - x\big| \lesssim |t-s|^{\frac{1}{p}}\big(1+\|\bfX\|\big)
$$
and 
$$
\sup_{x\in \bbR^d} \big|D^k \varphi_{ts}(x) - D^k \varphi_{s,s}(x)\big| \lesssim |t-s|^{\frac1p}e^{c_k |t-s|^\frac1p N},
$$
where $N = \Big[c\big(1+\|\bfX\|\big)^{-p}\Big]$.
\end{thm}

\medskip

\begin{Dem}
We work here with $\alpha=0$, so we know from the above computations that for $|t-s|^{\frac1p}\big(1+\|\bfX\|\big) \lesssim 1$, and all $0\leq k\leq (n+1)$, we have
\begin{equation}
\label{eq:bo}
\Big\| D^k\big(\mu^n_{ts}- \mu_{ts}\big)(x) \Big\|_\infty \lesssim |t-s|^{\frac{1+[p]}{p}}\big(1+\|\bfX\|\big)^{[p]+1}.
\end{equation}
This implies that for all $0\le k \le n$ the function $D^k_{t,s}$ is Lipschitz continuous, with Lipschitz constant not greater than a constant multiple of $|t-s|^{\frac{1+[p]}{p}}\big(1+\|\bfX\|\big)^{[p]+1}$. It follows from this fact and the proof of Theorems \ref{thm:local} and \ref{thm:main} that there exists some maps $A^k_{t,s}$, such that $D^k \mu^n_{ts}$ converges uniformly to $A^k_{t,s}$ as $n$ goes to $\infty$. One then needs to prove that the $A^k$ are indeed the $k$-th derivative of $\varphi$ and get the bounds of the statement. Note that the small time bounds are direct consequences of equation \eqref{eq:bo} and Remark \ref{rem:derivatives1} once we know that the $D^k\mu^n_{ts}$ converge.

\medskip

We have
$$
\mu^n_{ts}(x+h)-\sum_{0\le j\le  k} \frac{1}{j!}D^k \mu_{t,s}(x) \cdot h^{j} = \frac1{k!}\int_0^1 \dd \lambda D^{k+1}\mu_{ts}^n(\lambda h+  x)\cdot ((1-\lambda)h)^{k} h,
$$
where $h^j = \underbrace{(h,\cdots,h)}_{j \text{ times}}$. Hence, thanks to Remark \ref{rem:derivatives3}, for all $|t-s|^{\frac{1}{p}}\big(1+\|\bfX\|\big)\lesssim 1$, the maps $D^{k+1}\mu^n_{ts}$ are bounded, \emph{uniformly in $n$}, and
$$
\left|\mu^n_{t,s}(x+h)-\sum_{0\le j\le  k} \frac{1}{j!} D^k \mu^n_{t,s}(x) \cdot h^{j}\right| \lesssim |t-s|^{\frac{1+[p]}{p}} \big(1+\|\bfX\|\big)^{[p]+1} |h|^{k+1}.
$$
The previous bound allows us to send $n$ to $\infty$, and to get, as a consequence, that $A^k_{ts} = D^k \varphi_{ts}$. The construction of the global in time flow and its derivatives is done by gluing all these local flows, as above. 

\ssk

We now turn to the global bounds. As previously let $N$ be the least integer such that $T^{\frac1p}N^{-\frac1p}\big(1+\|\bfX\|\big) \lesssim 1$, where the implicit multiplicative constant is chosen such that all the previous bounds hold. Setting $t_i := \frac{i}{N}(t-s)+s$, one can use the local in time bounds on some time interval of length $(t_{i+1}-t_i)$. We have
\begin{align*}
\varphi_{t_i s}(x) - x =& \varphi_{t_i t_{i-1}}\big(\varphi_{t_{i-1} s}(x)\big)- \mu_{t_i t_{i-1}}\big(\varphi_{t_{i-1}s}(x)\big) \\
& +
\mu_{t_i t_{i-1}}\big(\varphi_{t_{i-1} s}(x)\big)  - \varphi_{t_{i-1} s}(x)\\
& + 
\varphi_{t_{i-1} s}(x)-x
\end{align*}
Hence, if one sets $R^0_i := \sup_{x} \big|\varphi_{t_i,s}(x)-x\big|$, one has
$$
R^0_{i} \le R^0_{i-1} + C|t-s|^{\frac{1}{p}} \lesssim i |t-s|^{\frac{1}{p}}.
$$
Similarly, we have
\begin{align*}
D\varphi_{t_i s}(x) - \textrm{Id} =& \Big(D\varphi_{t_i t_{i-1}}\big(\varphi_{t_{i-1} s}(x)\big)- \mu_{t_i t_{i-1}}\big(\varphi_{t_{i-1} s}(x)\big) \Big) D\varphi_{t_{i-1} s}(x) \\
& +
\Big(D\mu_{t_i t_{i-1}}\big(\varphi_{t_{i-1} s}(x)\big) - \textrm{Id}\Big) D\varphi_{t_{i-1} s}(x)\\
& + 
D\varphi_{t_{i-1} s}(x) - \textrm{Id}.
\end{align*}
Again,given the choice of $N$, one have can use all the local bounds on $\varphi$, $D\varphi$, and $\mu$ and $D\mu$, and setting $R^1_i := \sup_x \big|D\varphi_{t_i,s}(x) - \textrm{Id}\big|$, one has
$$
R^1_i \le C|t-s|^{\frac1p} + R^1_{i-1}\big(1+|t-s|\big)^{\frac{1}{p}}, 
$$
and 
$$
\sup_{x}\Big| D\varphi_{s,t}(x) - \textrm{Id}\Big| \lesssim |t-s|^{\frac1p}e^{c_1 |t-s|^\frac1p N}.
$$
One obtains the bounds for the higher order derivatives using Fa\`a di Bruno formula.
\end{Dem}

\bigskip

\subsection{Finite $p$-variation rough paths}
\label{SubsectionControls}

It is well-known the global bound for the differential of the flow, or the global bound for the flow for vector field with linear-growth, is not not good -- \cite{CCL}, \cite{riedel}, \cite{FrizHairer}. Indeed, in the setting of weak geometric Hölder $p$-rough paths, $N\sim \big(1+\|\bfX\|\big)^p$, and for a Gaussian rough path $\bfX$, the quantity $\|\bfX\|$ only has Gaussian tail and $\bbE\big[e^{c\|\bfX\|^p}\big] = +\infty$ for any $p>2$ and any positive constant $c$. To derive some moment bounds of solutions of rough differential equations, one need more advanced tools; we recall them here for the reader's convenience.

\medskip

\begin{definition*}
A weak geometric continuous rough path with finite $p$-variation is a continuous $[p]$-level weak geometric rough path such that 
\[\|\bfX\|_{[0,T],p-\mathrm{var}} := \sum_{i=1}^{[p]} \sup_{\pi \text{ partition of } [0,T]}\left(\sum_{(t_k,t_{k+1})\in \pi} |X^i_{t,s}|^{\frac{p}{i}}\right)^{\frac{1}{p}} <+\infty\]
Set 
\[
w(t,s) := \|\bfX\|_{[s,t],p-\mathrm{var}}^p.
\]
\end{definition*}

\medskip

If $\bfX$ is a weak geometric continuous rough path with finite $p$-variation then $w$ is a control; it is in particular increasing in its two variables, super-additive and continuous on the diagonal. Note also that a weak geometric Hölder $p$-rough path is always of finite $p$-variation since
\[
w(t,s) \le |t-s| \big(1+\|\bfX\|\big)^p.
\] 
The advantage of using the $p$-variation norm instead of the Hölder norm is related to integrability properties for random rough paths. 

\medskip

\begin{defn}\label{def:accumulation}
Given $\beta>0$ define $\tau^\beta_{0}=0$ and 
\[
\tau^\beta_{i+1} = \inf\{t \in [\tau^\beta_i,T] : w(\tau_i^\beta,t)\ge \beta\}\wedge T.
\]
The quantity $N_\beta := \sup\{ i \ge 0 : \tau^\beta_i < T\}$ is called the \textsf{\textbf{local accumulated variation of $\bfX$}}. 
\end{defn}

\medskip

The following result combines results from Friz and Victoir \cite{FrizVictoirBook} and Cass, Litterer and Lyons \cite{CCL}

\medskip

\begin{theorem*}\label{thm:CCL}
Let $\beta>0$, $p\ge 2$ and let $X$ be a centered Gaussian process defined over some finite interval $[0,T]$. Suppose that the covariance function is of finite $\rho$-two dimensional variation for some $\rho\in (1,2)$. Then for any $p\in(2\rho,4)$, $X$ can be lifted as a level-$[p]$ weakly geometric continuous finite $p$-variation rough path, and for $\beta>0$, the process $N_{\beta}^{\frac{1}{\rho}}$ has a Gaussian tails, namely there exists a constant $\mu>0$ such that
\[\bbE\left[ \exp\left(\mu N_\beta^{\frac{2}{\rho}}\right)\right] < +\infty.\]
In particular, for $p\in(2\rho,4)$, and for any constant $C>0$,
\[\bbE\left[\exp\left(C N_\beta\right)\right] \lesssim 1.\]
\end{theorem*}

\medskip

Friz and Riedel gave in \cite{riedel} what is now the classical proof of this result, based on Borell's isoperimetric inequality in Gaussian spaces. Cass and Ogrodnik \cite{CassOgrodnik} use heat kernel estimates as a substitute to isoperimetry to prove a similar result for Markovian rough paths. Compare the following definition to definition \ref{DefnSolutionFlow}.

\medskip

\begin{definition*}
A flow $\varphi : \Delta_T \times \RR^d \mapsto \RR^d$ is said to be a \textbf{\textsf{solution flow to the rough differential equation \eqref{EqRDEFlow}}} if there exists an exponent $\eta>1$ such that one can associate to any positive radius $R$ two positive constants $C_R$ and $\eps$, independent of $\bfX$, such that one has 
\begin{equation}
\label{EqDefnSolFlowVar}
\sup_{ x\in B(0,R)} \, \big|\varphi_{ts}(x) - \mu_{ts}(x) \big| \le C_{R}\, w(t,s)^{\eta},
\end{equation}
whenever $w(t,s)\leq \eps$.
\end{definition*}

\medskip

\begin{thm}\label{thm:main_variation} 
Let $\bfX$ be a weak geometric continuous rough path with finite $p$-variation. Let $V_0$ and $(V_1,\dots,V_\ell)$ satisfy \textbf{\textsf{Assumption 1}} and \textbf{\textsf{Assumption 2}}. There exists a unique global in time solution flow $\varphi$ to the rough differential equation \eqref{EqRDEFlow}.   \vspace{0.1cm}
\begin{itemize}
   \item One can choose $\eta=\frac{1+[p]}{p}$, $\eps = c_1 $ and $C_{R} = c_2 (1+R)^\alpha $, for some positive universal constants $c_1, c_2$, in the defining identity \eqref{EqDefnSolFlowVar}.   \vspace{0.1cm}
   
   \item One has for all $f\in C^{[p]+1}_b$ and all $w(t,s) \le \eps$ the estimate
\begin{multline*}
\sup_{x\in B(0,R)} \bigg|f\circ\varphi_{t,s}(x) - \Big\{
 f(x) + (t-s)V_{0}(s,\cdot)f + \sum_{k=1}^{[p]} \sum_{I \in \left\{0,\cdots, \ell\right\}^{k}} X^{k,I}_{t,s}V_{I}(s,\cdot)f\Big\}(x) \bigg|
 \\ \lesssim \|f\|_{C^{[p]+1}_b} (1+R)^{\alpha([p]+1)}w(t,s)^{\frac{[p]+1}{p}}.
\end{multline*}   
When $f = \textrm{Id}$, one can replace $(1+R)^{\alpha([p]+1)}$ by $(1+R)^\alpha$ and $\|f\|_{C^{n}_b}$ by $1$ in the previous bound.   \vspace{0.1cm}

   \item  The map that associates $\varphi$ to $\bfX$ is continuous from the set of weak geometric continuous rough paths with finite $p$-variation into the set of continuous flows endowed with the topology of uniform convergence on bounded sets.   \vspace{0.1cm}
   
   \item Finally, there exists $\beta>0$ and $c_3>0$ such that one has for all $(t,s)\in \Delta_T$,
   $$
   \underset{x\in B(0,R)}{\sup}\, \big|\varphi_{ts}(x)-x\big| \lesssim 
   \left\{\begin{aligned} 
	   &  \left( \left((1+R)^{1-\alpha}+c_4 w(t,s)^{\frac{1}{p}}N_{\beta}^{1-\frac1p}\right)^{\frac1{1-\alpha}}-(1+R)\right), & \text{if } \alpha<1\\
	   &(1+R)w(t,s)^{\frac{1}{p}}e^{c_3 N_\beta}, & \text{if } \alpha=1.
\end{aligned}
\right.
   $$
\end{itemize}
\end{thm}

\medskip

One gets back Theorem \ref{thm:main} when $\bfX$ is an Hölder $p$-rough path, with $N$ replaced by $N_\beta$.

\medskip

\begin{Dem}
The proof follows exactly the same steps as the proofs of Theorems \ref{thm:local} and Theorem \ref{thm:main}. We give here the main changes and leave the computations to the reader.

\medskip 

First, there is no loss of generality in assuming that $|t-s| \le w(t,s)$; replace if necessary $w(t,s)$ by $\tilde w(t,s) = |t-s|+ w(t,s)$. Set 
\[
C(t,s,\bfX) := \sum_{k=1}^{[p]} w(t,s)^{\frac{k}{p}}.
\]
One can replace the constant $C_{|t-s|,\|\bfX\|}$ by $C(t,s,\bfX)$ in Lemma \ref{lemma:gronwall_bounds} and Remark \ref{rem:derivatives1} as soon as $w(t,s) \le 1$; this ensures that $C(t,s,\bfX) \lesssim 1$. Lemma \ref{lemma:taylor} remains the same as it relies only on algebraic manipulations. In Lemma \ref{lemma:bound_remainder}, one has to assume that $w(t,s,\bfX) \le 1$, and one can replace $\big(1+\|\bfX\|\big)^{[p]+1}|t-s|^{\frac{1+[p]}{p}}$ in the estimates by $w(t,s)^{\frac{1+[p]}{p}}$ -- recall $|t-s|\le w(t,s)$. The same replacement is done in Corollary \ref{corollary:taylor} and Remark \ref{rem:derivatives2}. Finally, using the inequality $|t-s|\le w(t,s)$ and the fact that the real-valued functions $u\to w(t,u)$ and  $u\to w(u,s)$ are increasing, one can also replace $\big(1+\|\bfX\|\big)^{[p]+1}|t-s|^{\frac{1+[p]}{p}}$ by $w(t,s)^{\frac{1+[p]}{p}}$ in Proposition \ref{prop:sew}.

\smallskip

The proofs of Proposition \ref{prop:sewing1} and Proposition \ref{prop:sewing2} are a bit different, but the spirit is the same. The main difference is that one cannot say immediately that $w\big(t,\frac{t+s}2\big) \le \frac12 w(t,s)$. But given $(t,s) \in \Delta_T$, there exists $\tilde u\in (s,t)$ such that $w(t,u) = w(u,s) \le \frac12 w(t,s)$. Consider any sequence of embedded partitions $(\pi^n)_{n\in \bbN}=\Big(\big(t^n_{i}\big)_{i\in \{0,\cdots,n\}}\Big)_{n\in\bbN}$ with mesh going to $0$. One proves by induction the existence of constants $0<\beta\le 1$ and $L>0$ such that for $w(t,s) \le \beta$, one has for all $k\le n$, 
\[
\sup_{x \in B(0,R)}\,\Big|\mu_{t^k_k,t^k_{k-1}}\circ \cdots \circ \mu_{t^k_{1},t^k_0}(x) - \mu_{t,s}(x)\Big|\le L\, (1+R)^\alpha \, w(t,s)^{\frac{[p]+1}{p}} 
\]
Let the integer $0\leq i_0\leq n$ be such that $t^{n+1}_{i_0} \le \tilde u < t^{n+1}_{i_0+1}$. One closes the induction and proves the following bound for all $n\in \bbN$ by taking $u = t^{n+1}_{i_0+1}$, using the fact that 
\[
w(t^{n+1}_{i_0+1},\tilde u) + w(\tilde u,t^{n+1}_{i_0}) \underset{n\to \infty}{\longrightarrow} 0.
\]
The same trick holds for the proof of Proposition \ref{prop:sewing2}, assuming that 
\[
w(t,s)(1+R)^{\frac{\alpha}{1+[p]}} \le \beta.
\]
One can again replace in Proposition \ref{prop:sewing1},Propositino \ref{prop:sewing2} and Remark \ref{rem:derivatives3} $\big(1+\|\bfX\|\big)^{[p]+1}|t-s|^{\frac{1+[p]}{p}}$ by $w(t,s)^{\frac{1+[p]}{p}}$.

\medskip

For the proof of the local existence, one can proceed as in Lemma 2.1 of \cite{DeyHofGubTin}, and as in the proof of Theorem \ref{thm:local}. Let $\big((t^n_i)_{i\in \{0,\cdots,n\}}\big)_{n\in\bbN}$ be the sequence of dyadic partitions. Remark that since $w$ is superadditive, there exists $i$ such that $w(t^n_{i+1},t^n_{i-1}) \le (2^{n}-1)^{-1}w(t,s)$. Define the partition $\widehat \pi= \big\{s=t_0, < \cdots < t^n_{i-1} < t^n_{i+1} <\cdots < t^{n}_{2^n} = t\big\}$and set $M^n_{t,s} := \mu^{n}_{t,s} - \mu_{t,s}$, and 
\[
\widehat M^n_{t,s} := \mu_{t^{n}_{2^n},t^{n}_{2^n-1}}\circ\cdots \circ \mu_{t^n_{i+2},t^n_{i+1}}\circ \mu_{t^n_{i+1},t^n_{i-1}}\circ \mu_{t^n_{i-1},t^n_{i-2}}\circ \mu_{t^n_1, t^n_0} - \mu_{t,s}.
\] 
We have
\begin{multline*}
\widehat M_{s,t} -M^n_{t,s}
= 
\Big\{\mu_{t^{n}_{2^n},t^{n}_{2^n-1}}
	\circ
\cdots 
	\circ 
\mu_{t^n_{i+2},t^n_{i+1}} 
	\circ (\mu_{t^n_{i+1},t^n_i}\circ \mu_{t^n_{i},t^n_{i-1}})\\ 
	- \mu_{t^{n}_{2^n},t^{n}_{2^n-1}}
	\circ
\cdots 
	\circ 
\mu_{t^n_{i+2},t^n_{i+1}}\circ (\mu_{t^n_{i+1},t^n_{i-1}} )\Big\} \circ\mu_{t^n_{i-1},t^n_{i-2}}\circ \mu_{t^n_1, t^n_0}(x).
\end{multline*}
The induction hypothesis and the bound $w(t,s)(1+R)^{\frac\alpha{1+[p]}} \le \beta$, then give
\[
\big|\widehat M^n_{s,t} -M^n_{t,s}\big| \lesssim (2^{n}-1)^{-\frac{[p]+1}{p}}(1+R)^\alpha w(s,t)^{\frac{[p]+1}{p}}.
\]
Repeating this operation until we get the trivial partition of $[s,t]$ we see that
\[
M^n_{t,s} = \sum_{k=0}^{2^n} \rho^k_{t,s},
\]
with 
$$
\big|\rho^k_{t,s}(x)\big| \lesssim (1+R)^\alpha \, w(t,s)^{\frac{[p]+1}{p}} \,(2^{n}-k)^{-\frac{[p]+1}{p}}.
$$
Here we crucially use the fact that the composition of the flows are globally Lipschitz continuous, uniformly in $n$. Hence $M^n$ converges uniformly to a limit $\varphi_{t,s } - \mu_{t,s}$ and
\[
\sup_{x\in B(0,R)} |\varphi_{t,s}(x)-\mu_{t,s}(x)| \lesssim (1+R)^\alpha \sum_{i\ge 0} i^{-\frac{[p]+1}{p}} w(t,s)^{\frac{[p]+1}{p}} \lesssim (1+R)^\alpha w(t,s)^{\frac{[p]+1}{p}}.
\]
The remainder of the proof follows easily from the proof of Theorem \ref{thm:local} and Theorem \ref{thm:main}. Indeed, by construction, $\varphi$ is a flow for all dyadic points, and then by continuity for all points, and thanks to the continuity of $\mu$ with respect to $\bfX$, $\varphi$ is continuous with respect to $\bfX$. 

\smallskip

Note also that thanks to the superadditivity property of the control, one has
\[\sum_{k=0}^{2^n-1} w(t^n_{i+1},t^n_i)^{\frac{[p]+1}{p}} \lesssim \max_{i\in\{0,\cdots,2^n-1\}} w(t^n_{i+1},t^n_i)^{\frac{[p]+1-p}{p}} w(t,s),
\]
and since $w$ is continuous on the diagonal, the above sum goes to $0$ as $n$ goes to infinity. Local uniqueness of the flow  follows -- see Equation \eqref{eq:local_uniqueness}. 

\smallskip

The proof of global existence is similar to the proof of Theorem \eqref{thm:main}. Use the sequence of times $(\tau^\beta_i)_{i\in \bbN}$ from definition \ref{def:accumulation}. We have
\begin{align*}
\varphi_{\tau^\beta_i,s}(x) -x &= \varphi_{\tau^{\beta}_{i}\tau^{\beta}_{i-1}}\big(\varphi_{\tau^{\beta}_{i-1}s}(x)\big) - \mu_{\tau^{\beta}_{i}\tau^{\beta}_{i-1}}\big(\varphi_{\tau^\beta_{i-1} s}(x)\big)   \\
&\quad+\mu_{\tau^{\beta}_{i}\tau^{\beta}_{i-1}}\big(\varphi_{\tau^\beta_{i-1} s}(x)\big) - \varphi_{\tau^\beta_{i-1} s}(x)   \\
&\quad+\varphi_{\tau^\beta_{i-1} s}(x) - x.
\end{align*}
Define $R_i := \sup_{x\in B(0,R)} \big|\varphi_{\tau^\beta_{i} s}(x)  - x \big|$ and $R_0 = 0$. The fourth item of the statement follows then from the induction relation
\[R_i \le R_{i-1} + w(\tau^{\beta}_i,\tau^\beta_{i-1})^{\frac{[p]+1}{p}}(1+R+R_i)^{\alpha} + C(\tau^\beta_{i},\tau^\beta_{i-1},\bfX)(1+R+R_{i-1})^\alpha.\]
Since $w(\tau_i^\beta, \tau^\beta_{i-1}) \le \beta $, one has $C(\tau^\beta_i,\tau^\beta_{i-1},\bX) \lesssim w(\tau^\beta_i,\tau^\beta_{i-1})^{\frac1p}.$, hence
\[R_{i}-R_{i-1} \lesssim (1+R+R_{i-1})^{\alpha} \big(w(\tau^\beta_i,\tau^\beta_{i-1})^{\frac{[p]+1}p}+w(\tau^\beta_i,\tau^\beta_{i-1})^{\frac1p}\big).\]
When $\alpha=1$ one end up with the following bound :
\[R_{i} \lesssim  R_{i-1} + (1+R)w(t,s)^{\frac1p}.\]
When $\alpha<1$, one ends up with 
\[R_N \lesssim \bigg(\Big((1+R)^{1-\alpha} + \frac{1}{1-\alpha}\sum_{i=1}^N\big(w(\tau_i,\tau_{i-1})^{\frac{[p]+1}{p}}+w(\tau_i,\tau_{i-1})^{\frac{1}{p}}\big) \Big)^{\frac1{1-\alpha}} - (1+R)\bigg).\]
By using Jensen formula, one finally has the bound
\[\sum_{i=1}^N\big(w(\tau_i,\tau_{i-1})^{\frac{[p]+1}{p}}+w(\tau_i,\tau_{i-1})^{\frac{1}{p}}\big) \lesssim N^{1-\frac1p} w(t,s)^{\frac1p},\]
which ends the proof.
\end{Dem}

\medskip

\begin{thm}\label{thm:derivatives_variation}
Let $p>2$ and $\bfX$ be a weak geometric continuous finite $p$-variation rough path and let $(V_0,\cdots V_\ell)$ which satisfies Assumption \ref{hyp:derivatives}.
Let $\varphi$ stands for the solution flow to the rough differential equation \eqref{EqRDEDerivativeFlow}. Then each $\varphi_{t,s}$ is of class $C^n$, has linear growth and bounded derivatives. Furthermore for a suitable positive constant $\epsilon_3$, independent of $\bfX$, and $w(t,s) \leq \epsilon_3$, we have
$$
\Big\| D^k\varphi_{ts} -D^k \mu_{ts}\Big\|_\infty \lesssim w(t,s)^{\frac{1+[p]}{p}},
$$
for all $0\leq k\leq n$. Finally there exists $\beta>0$ and some positive constants $c_1,\cdots,c_n$, independent of $\bfX$, such that \emph{for all} $0\le s \le t \le T$, and every $1\leq k\leq n$, we have
$$
\sup_{x\in \bbR^d} \big| \varphi_{ts}(x) - x\big| \lesssim w(t,s)^{\frac{1}{p}}N_\beta
$$
and 
$$
\sup_{x\in \bbR^d} \big|D^k \varphi_{ts}(x)-D^k \varphi_{s,s}(x)\big| \lesssim w(t,s)^{\frac{1}{p}}e^{c_k  N_\beta}.
$$
\end{thm}

\medskip

\begin{Dem}
We refer to the proof of Theorem \ref{thm:derivatives}. The first bound of the theorem is a direct application of Theorem \ref{thm:main_variation} with $\alpha=0$. For the existence of derivatives and the associated bounds, one can mimic the proof of Theorem \ref{thm:derivatives} by replacing $|t-s|^{\frac1p}\big(1+\|\bfX\|\big)^p$ by $w(t,s)$. The proof of the global bound is done in the same way as the proof of Theorem \ref{thm:main_variation}.
\end{Dem}

\bigskip

\section*{Acknowledgements}
Most of this work was done while R. Catellier was working at University of Rennes 1 (IRMAR -- UMR 6625), with the support of the Labex Centre Henri Lebesgue.

\bigskip
\bigskip

\noindent \textcolor{gray}{$\bullet$} {\sf I. Bailleul} - {\small Institut de Recherche Mathematiques de Rennes, 263 Avenue du General Leclerc, 35042 Rennes, France.} {\it ismael.bailleul@univ-rennes1.fr}   \vspace{0.3cm}

\noindent \textcolor{gray}{$\bullet$} {\sf R. Catellier} - {\small Universit\'e C\^otes d'Azur, LJAD, Nice, France.} {\it remi.catellier@unice.fr}


\bigskip
\bigskip



\begin{thebibliography}{AAA}

\bibitem{LNRoughPaths}
Bailleul, I.,
\newblock A flows-based approach to rough differential equations.
\newblock https://perso.univ-rennes1.fr/ismael.bailleul/files/M2Course.pdf.

\bibitem{BailleulSeminaire}
Bailleul, I.
\newblock Flows driven by Banach space walued rough paths.
\newblock {\em S\'eminaire de Probabilit\'es}, Vol. {XLVI}:195--205, 2015.

\bibitem{BailleulRMI}
Bailleul, I.
\newblock Flows driven by rough paths.
\newblock {\em Revista Mat. Iberoamericana}, {31}(3):901--934, 2015.

\bibitem{BailleulCatellierMeanField}
Bailleul, I. and Catellier, R.,
\newblock Limit theorems for systems of mean field rough differential equations.
\newblock Preprint, 2018.

\bibitem{BaudoinBookEMS}
Baudoin, F.,
\newblock Diffusion Processes and Stochastic Calculus.
\newblock {\em EMS Textbooks in Mathematics}, 2014.

\bibitem{BoutaibGyurkoLyonsYang}
Boutaib, Y. and Gyurko, L.G. and Lyons, T. and Yang, D.,
\newblock Dimension-free Euler estimates of rough differential equations.
\newblock {\em Rev. Roumaine Math. Pures Appl.}, 2014.

\bibitem{LyonsStFlour} 
Caruana, M. and L\'evy, T. and Lyons, T., 
\newblock Differential equations driven by rough paths.
\newblock {\em Lect. Notes in Math.}, {\bf 1908}, 2004.

\bibitem{CCL}
Cass, T.,  Litterer, C. and Lyons, T.,
\newblock
Integrability and tail estimates for Gaussian rough differential equations
\newblock
{\em The Annals of Probability} {\bf 41}(4):3026--3050, 2013. 

\bibitem{CassOgrodnik}
Cass, T. and Ogrodnik, M.
\newblock
Tail estimates for Markovian rough paths
\newblock
{\em The Annales of Probability}, {\bf 45}(4):2477--2504,.

\bibitem{CassWeidner}
Cass, T. and Weidner, M.,
\newblock Tree algebras over topological vector spaces in rough path theory
\newblock arXiv:1604.07352v2, 2016.

\bibitem{LejayCoutinOmega}
Coutin, L. and Lejay, A. 
\newblock Sensitivity of rough differential equations: an approach through the Omega lemma.
\newblock {\em J. Diff. Eq.}, {\bf 264}(6):3899--3917, 2018.

\bibitem{CoutinLejayPerturbed}
L. Coutin and A. Lejay, 
\newblock Perturbed linear rough differential equations,
\newblock {\em Ann. Math. Blaise Pascal}, {\bf 21}(1):103--150, 2014.

\bibitem{Davie}
Davie, A.M.,
\newblock Differential equations driven by rough paths: an approach via discrete approximations.
\newblock {\em Appl. Math. Res. Express}, 9--40, 2007.

\bibitem{DeyHofGubTin}
A. Deya, M. Gubinelli, M. Hofmanova and S. Tindel
\newblock
A priori estimates for rough PDEs with application to rough conservation laws.
\newblock arXiv:1604.00437, 2016.


\bibitem{FrizHairer}
Friz, P. and Hairer, M.,
\newblock
 A Course on Rough Paths,
 \newblock
 {\em Spinger, Universitext}, 2014. 
 
\bibitem{riedel}
Friz, P. and Riedel, S. 
\newblock
Integrability of (non-)linear rough differential equations and integrals. 
\newblock{\em Stochastic Analysis and Applications}, {\bf 31}(2):336--358, 2013.

\bibitem{FrizVictoirBook}
Friz, P. and Victoir, N.,
\newblock Multidimensional stochastic processes as rough paths.
\newblock {\em Cambridge studies in advanced Mathematics}{\bf 120}, 2010.

\bibitem{LejayControl}
Lejay, A.,
\newblock Controlled differential equations as Young integrals: a simple approach.
\newblock {\em J. Differential Equations}, {\bf 249}:1777--1798, 2010.

\bibitem{LejayGlobal}
Lejay, A.,
\newblock Global solutions to rough differential equations with unbounded vector fields.
\newblock {\em S\'em. Probailit\'es}, Vol. {\bf XLIV}:215--246, 2012.

\bibitem{Lejay}
Lejay, A.,
\newblock On rough differential equations.
\newblock {\em Elec. J. Probab.}{\bf 14}(12), 341--364, 2009.

\bibitem{Lyons98}
Lyons, T.,
\newblock Differential equations driven by rough signals.
\newblock {\em Revista Mat. Iberoamaricana}, {\bf 14}(2):215--310, 1998. 

\bibitem{LyonsQian}
Lyons, T. and Qian, Z.,
\newblock System control and rough paths.
\newblock {\em Oxford Mathematical Monographs}, 2002.

\bibitem{LyonsYang}
Y. Boutaib and Lyons, T. and Yang, D.,
\newblock Dimension-free Euler estimates of rough differential equations.
\newblock {\em Rev. Roumanine Math. Pures Appl.}, {\bf 59}(1):25--53, 2014.

\end{thebibliography}
\end{document}